%% file: FlexSummary_V15_2016_Rev2-Final.tex
\documentclass[opre,nonblindrev]{informs3aa} %

\OneAndAHalfSpacedXI 

\usepackage{mathtools}
\usepackage{amssymb}
\usepackage{comment}
\usepackage{subfigure}
\usepackage{bm}
\usepackage{url}
\usepackage{graphicx}
\usepackage{epstopdf}
\usepackage{epsfig,psfrag}
\usepackage{verbatim} 
\usepackage{subfigure}
\usepackage{algorithmic}
\usepackage{enumerate}
\usepackage{color}
\usepackage{pdfsync}
\usepackage{hyperref}
\usepackage{array}

\theoremstyle{plain}
\newtheorem{theorem}{Theorem}[section]
\newtheorem{thm}[theorem]{Theorem}	
\newtheorem{lemma}[theorem]{Lemma}

\newtheorem{prop}[theorem]{Proposition}
\newtheorem{proposition}[theorem]{Proposition}

\newtheorem{definition}[theorem]{Definition}
\newtheorem{defn}[theorem]{Definition}

\newtheorem{conjecture}[theorem]{Conjecture}

\newtheorem{cond}[theorem]{Condition}

\newtheorem{assumption}[theorem]{Assumption}

\include{macroDef}

\usepackage[square]{natbib}
 \bibpunct[, ]{[}{]}{,}{a}{}{,}%
 \def\bibfont{\small}%
 \def\bibsep{\smallskipamount}%
 \def\bibhang{24pt}%
\def\ell{l}
\def\s{s}

\def\bdel#1{}

\def\nln{\nonumber \\}

\begin{document}


\RUNAUTHOR{Tsitsiklis and Xu}

\RUNTITLE{Flexible Queueing Architectures}

\TITLE{Flexible Queueing Architectures}

\ARTICLEAUTHORS{
\AUTHOR{John N.~Tsitsiklis}
\AFF{LIDS, Massachusetts Institute of Technology, Cambridge, MA 02139, \EMAIL{jnt@mit.edu}} 
\AUTHOR{Kuang Xu}
\AFF{Graduate School of Business, Stanford University, 
Stanford, CA 94305, \EMAIL{kuangxu@stanford.edu}} 
} 

%
%
%
%
%

\ABSTRACT{
We study a multi-server model with $n$ \emph{flexible} servers and $n$ queues, connected through a bipartite graph, where the level of flexibility is captured by an upper bound on the graph's average degree, $d_n$. Applications in content replication in data centers, skill-based routing in call centers, and flexible supply chains are among our main motivations.  

We focus on the scaling regime where the system size $n$ tends to infinity, while the overall traffic intensity stays fixed. We show that a large capacity region and an asymptotically vanishing queueing delay  are simultaneously achievable even under limited flexibility ($d_n \ll n$). Our main results demonstrate that, when $d_n\gg \ln n$, a family of expander-graph-based flexibility architectures {has} a capacity region that is within a constant factor of the maximum possible, while simultaneously ensuring a diminishing queueing delay for \emph{all} arrival rate vectors in the capacity region. 
Our analysis is centered around a new class of virtual-queue-based scheduling policies that rely on dynamically constructed job-to-server assignments on the connectivity graph. 
For comparison, we also analyze a natural family of modular architectures, which is simpler but has provably weaker performance. \footnote{May 2015; revised October 2016. A preliminary version of this paper appeared at Sigmetrics 2013, \cite{TX13SIG}; the performance of the architectures proposed in the current paper is significantly better than the one in \cite{TX13SIG}. This research was supported in part by the NSF under grant CMMI-1234062.}}


\maketitle

\section{Introduction}
\label{sec:randomGraphintro}
At the heart of a number of modern queueing networks lies the problem of allocating processing resources (manufacturing plants, web servers, or call-center staff) to meet multiple types of demands that arrive dynamically over time (orders, data queries, or customer inquiries). It is usually the case that a \emph{fully flexible} or \emph{completely resource-pooled} system, where every unit of processing resource is capable of serving all types of demands, delivers the best possible performance. Our inquiry is, however, motivated by the unfortunate reality that such full flexibility is often infeasible due to overwhelming implementation costs (in the case of a data center) or human skill limitations (in the case of a skill-based call center).

What are the key benefits of flexibility and resource pooling in such queueing networks? Can we harness the same benefits even when the degree of flexibility is \emph{limited}, and how should the network be designed and operated? These are the main questions that we wish to address. While these questions can be approached from a few different angles, we will focus on the metrics of \emph{capacity region} and \emph{expected queueing delay}; the former measures the system's \emph{robustness}  against \emph{demand uncertainties}, i.e., when the arrival rates for different demand types are unknown or likely to fluctuate over time, while the latter is a direct reflection of \emph{performance}. Our main message is positive: in the regime where the system size is large, improvements in both the capacity region and delay are \emph{jointly achievable} even under very limited flexibility, given a proper choice of the architecture (interconnection topology) and scheduling policy. 

\begin{figure}[h]
\centering
\includegraphics[scale=.7]{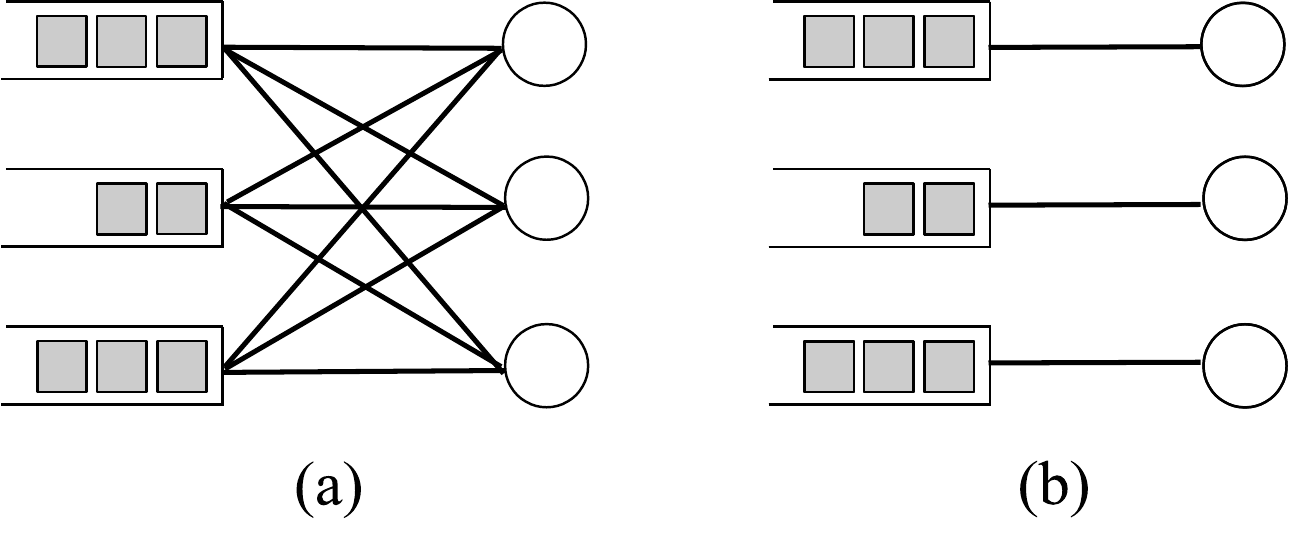}
\caption{Extreme cases of flexibility: $\ldn=n$ versus $\ldn=1$.}
\label{fig:2}
\vspace{-15pt}
\end{figure}

\bfpara{Benefits of Full Flexibility.} We begin by illustrating the benefits of flexibility and resource pooling {in a very simple setting.} Consider a system of $n$ servers, each running at rate $1$, and $n$ queues, where each queue stores jobs of a particular demand type. For each $ i \in \left\{1,\ldots, n\right\}$, queue $i$ receives an independent Poisson arrival stream of rate $\lambda_i$. The average arrival rate $\frac{1}{n}\sum_{i=1}^n \lambda_i$ is denoted by $\rho$, and is referred to as the \emph{traffic intensity}. The sizes of all jobs are independent and exponentially distributed with mean $1$. 

For the remainder of this paper, we will use a \emph{measure of flexibility} given by the average number of servers that a demand type can receive service from, denoted by $\ldn$. Let us consider the two extreme cases: a \emph{fully flexible} system, with $\ldn=n$ (Figure \ref{fig:2}(a)), and an \emph{inflexible} system, with $\ldn=1$ (Figure \ref{fig:2}(b)). Fixing the traffic intensity $\rho<1$, and letting the system size, $n$, tend to infinity, we observe the following qualitative benefits of full flexibility:
 
\bfpara{1.~Large Capacity Region.}  In the fully flexible case and under any work-conserving scheduling policy\footnote{A work-conserving policy mandates that a server be always busy whenever there is at least one job in {some queue} to which it is connected.}, the \emph{collection of all jobs} in the system evolves as an $M/M/n$ queue, with arrival rate $\sum_{i=1}^n \lambda_i$ and service rate $n$. It is easy to see that the system is stable for all arrival rates that satisfy $\sum_{i=1}^n \lambda_i < n$. {In contrast,} in the inflexible system, since all $M/M/1$ queues operate independently, we must have $\lambda_i< 1$, for all $i$, in order to achieve stability. Comparing the two, we see that the fully flexible system {has} a much larger capacity region, and is hence more robust to uncertainties or changes in the arrival rates. 

\bfpara{2.~Diminishing Delay.} Let $W$ be the steady-state expected waiting time in queue (time from entering the queue to the initiation of service). As mentioned earlier, the total number jobs in the system for the fully flexible case evolves as an $M/M/n$ queue with traffic intensity $\rho<1$. It is not difficult to verify that for any fixed value of $\rho$, the expected total number of jobs in the queues is \emph{bounded above} by a constant independent of $n$, and hence the expected waiting time in queue satisfies $\EX{W} \to 0$, as $n\to \infty$.\footnote{The {fact that the expected waiting time vanishes asymptotically} follows from the bounded expected total number of jobs in steady-state, the {assumption} that the total arrival rate is $\rho n$, which goes to infinity as $n\to \infty$, and Little's Law.} In contrast, the inflexible system is simply a collection of $n$ independent $M/M/1$ queues, and hence the expected waiting time is $\EX{W}=\frac{\rho}{1-\rho}>0$, for all $n$. Thus, the expected delay  in {the} fully flexible system {\emph{vanishes asymptotically} as the system size increases}, but stays bounded away from zero in {the} inflexible system. 

\bfpara{Preview of Main Results.} Will the above benefits of fully flexible systems continue to be present if the system only has limited flexibiltiy, that is, if $\ldn \ll n$? The main results of this paper show that a large capacity region and {an asymptotically vanishing} delay  can still be \emph{simultaneously achieved}, even when $\ldn\ll n$. However, when flexibility is limited, the architecture and scheduling policy need be chosen with care. {We show that, when $d_n\gg \ln n$,  a family of expander-graph-based flexibility architectures {has the largest possible} capacity region, up to a constant factor, while simultaneously ensuring a diminishing queueing delay, of order $\ln n/d_n$ as $n \to \infty$, for \emph{all} arrival rate vectors in the capacity region (Theorem \ref{thm:expanderArch}). For comparison, we also analyze a natural family of modular architectures, which is simpler but has provably weaker performance (Theorems \ref{thm:ModCapacity} and~\ref{thm:randModular}).}

\subsection{Motivating Applications}

\begin{figure}[t]
\centering
\includegraphics[scale=.85]{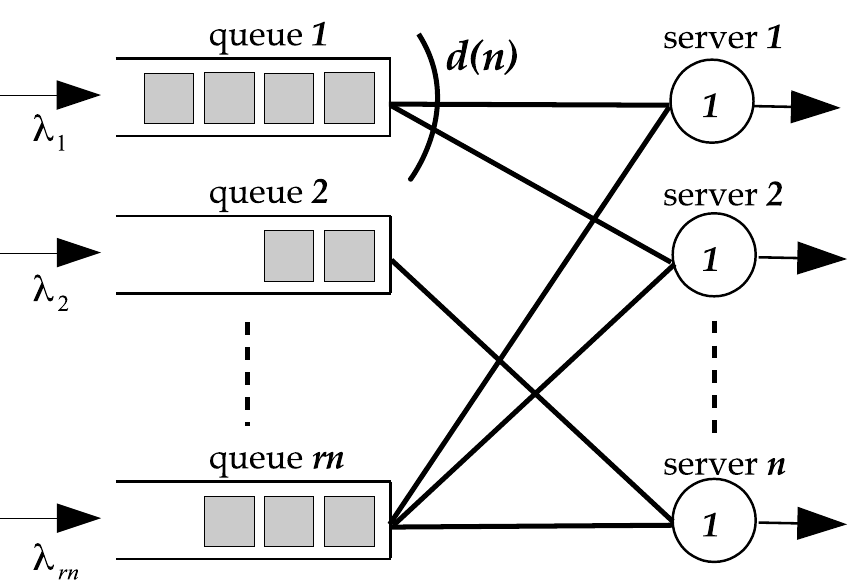}
\caption{A processing network with $rn$ queues and $n$ servers. 
}
\label{fig:random1}
\vspace{-15pt}
\end{figure}

We describe here several motivating applications for our model; Figure \ref{fig:random1} illustrates the overall architecture that they share. {\bf Content replication} is commonly used in data centers for bandwidth intensive operations such as database queries \cite{SAG06} or video streaming \cite{LLM12}, by hosting the same piece of content on multiple servers. Here, a server corresponds to a physical machine in the data center, and each queue stores incoming demands for a particular piece of content (e.g., a video clip). A server $j$ is connected to queue $i$ if there is a copy of content $i$ on server $j$, and $\ldn$ {reflects} the average number of replicas per piece of content across the network. Similar structures also arise in {\bf skill-based routing in call centers}, where agents (servers) are assigned to answer calls from different categories (queues) based on their domains of expertise \cite{WW05}, and in {\bf process-flexible supply chains} \cite{JG95, SW11, CCTZ10b, IVS05, GB03}, where each plant (server) is capable of producing multiple product types (queues). In many of these applications, demand rates can be unpredictable and may change significantly over time; for instance, unexpected ``spikes'' in demand traffic are common in modern data centers \cite{KSG09}. These demand uncertainties make \emph{robustness} an important criterion for system design. These practical concerns have been our primary motivation for studying the {interplay} between robustness, performance, and the level of flexibility. 

\subsection{Related Research}

Bipartite graphs provide a natural model for capturing the relationships between demand types and service resources. It is well known in the supply chain literature that limited flexibility, corresponding to a sparse bipartite graph, can be surprisingly effective in resource allocation even when compared to a fully flexible system \cite{JG95,GB03, IVS05, CCTZ10b, SW11}. The use of sparse random graphs or expanders as flexibility structures to improve robustness has recently been studied in \cite{CTZ10a,chen2014optimal} in the context of supply chains, and in \cite{LLM12} for content replication. Similar to the robustness results reported in this paper, {these} works show that random graphs or expanders can accommodate a large set of demand rates. However, in contrast to our work, nearly all analytical results in this literature focus on static allocation problems, where one tries to match supply with demand in a single shot, as opposed to our model, where resource allocation decisions need to be made dynamically over time.

In the queueing theory literature, the models that we consider fall under the umbrella of multi-class multi-server systems, where a set of servers are connected to a set of queues through a bipartite graph. Under these (and similar) settings, complete resource pooling (full flexibility) is known to improve system performance \cite{MR98,HL99,BW01}. However, much less is known when only limited flexibility is available: systems with a non-trivial connectivity graph are extremely difficult to analyze, even under seemingly simple scheduling policies (e.g, first-come first-serve) \cite{TW08,VADW12}. Simulations in \cite{WW05} show empirically that limited cross-training can be highly effective in a large call center under a skill-based routing algorithm. Using a very different set of modeling assumptions, \cite{BRV11} proposes a specific chaining structure with limited flexibility, which is shown to perform well under heavy traffic. Closer to the spirit of the current work is \cite{TX12}, which studies a partially flexible system where a fraction $p>0$ of all processing resources are fully flexible, while the remaining fraction, $1-p$, is dedicated to specific demand types, and which shows an exponential improvement in delay scaling under heavy-traffic. However, both \cite{BRV11} and \cite{TX12} focus on the heavy-traffic regime, which is different from the current setting where traffic intensity is assumed to be fixed, and the analytical results in both works apply only to uniform demand rates. Furthermore, with a constant fraction of the resources being fully flexible,
 the average degree in \cite{TX12} scales linearly with the system size $n$, whereas here we are interested in the case of a much slower (sub-linear) degree scaling. 

At a higher level, our work is focused on the interplay between robustness, delay, and the degree of flexibility in a queueing network, which is much less studied in the existing literature, and especially for networks with a non-trivial interconnection topology.

On the technical end, we build on several existing ideas. The techniques of batching (cf.\ \cite{NMC07,TS08}) and the use of virtual queues (cf.\ \cite{MMA99,KS01}) have appeared in many contexts in queueing theory, but the specific models considered in the literature bear little resemblance to ours. The study of expander graphs has become a rich field in mathematics (cf.~\cite{hoory2006expander}), {but} we will refrain from providing a thorough review because only some elementary and standard properties of expander graphs are used in the current paper. 

We finally note that preliminary (and weaker) versions of some of the results were included in the conference paper \cite{TX13SIG}. 

\bfpara{Organization of the Paper.} We describe the model in Section \ref{sec:modNotRandomGraph}, along with the notation to be used throughout. The main {results} are provided in Section \ref{sec:mResults}. The construction and the analysis associated with the Expander architecture will be presented separately, in Section \ref{sec:VirtualQueues}.  We conclude the paper in Section \ref{sec:discFut} with a further discussion of the results as well as directions for future research. 

\section{Model and Metrics}
\label{sec:modNotRandomGraph}

\subsection{Queueing Model and Interconnection Topologies}

{\bf The Model.}  We consider a sequence of systems operating in \emph{continuous time}, indexed by the integer $n$, where the $n$th system consists of  {$rn$} queues and $n$ servers (Figure \ref{fig:random1}), {and where  $r$ is a constant that is held fixed as $n$ varies. For simplicity, we will set $r$ to 1 but} {note} that {all results} and arguments in this paper can be extended to the case of 
{general $r$} 
without  difficulty. 

{A}  \emph{flexible architecture} is represented by an $n \times n$ undirected bipartite graph $g_n =\left(E, I\cup J\right)$, where $I$ and $J$ represent the sets of queues and servers, respectively, and $E$ the set of edges between them.\footnote{For simplicity of notation, we omit the dependence of $E$, $I$, and $J$ on $n$.} We will also refer to $I$ and $J$ as the sets of left and right {nodes}, respectively. A server $j \in J$ is \emph{capable} of serving a queue $i \in I$, if and only if $\left(i,j\right) \in E$. We will use the following notation. 
\begin{enumerate}
\item Let $\mcal{G}_n$ be the set of all $n \times n$ bipartite graphs.
\item For $g_n \in \gcal_n$, let $\degr(g_n)$ be the average degree among the $n$ left {nodes, which is} the same as the average degree of the right {nodes}.
\item For a subset of {nodes}, $M \subset I\cup J$, let $g|_{M}$ be the graph induced by $g$ on the {nodes} in $M$. 
\item Denote by $\neig{i}$ the set of servers in $J$ connected to queue $i$, and similarly, by $\neig{j}$ the set of queues in $I$ connected to server~$j$.
\end{enumerate}

Each queue $i$ receives a stream of incoming jobs according to a Poisson process of rate $\lambda_{n,i}$, independent of all other streams, and we define $\blam_n = \left(\lambda_{n,1},\lambda_{n,2},\ldots, \lambda_{n,\llfl n\rrfl}\right)$, which is the {\bf arrival rate vector}. When the value of $n$ is clear from the context, we sometimes suppress the subscript $n$ and write $\blam =(\lambda_1,\ldots, \lambda_n)$ instead. The sizes of the jobs are exponentially distributed with mean $1$, independent from each other and from the arrival processes. All servers are assumed to be running at a constant rate of $1$. The system is assumed to be empty at time $t=0$.


Jobs arriving at queue $i$ can be assigned (immediately, or in the future) to an idle server $j \in \neig{i}$ to receive service. The assignment is \emph{binding}: once the assignment is made, the job cannot be transferred to, or simultaneously receive service from, any other server. Moreover, service is \emph{non-preemptive}:  once service is initiated for a job, the assigned server has to dedicate its full capacity to this job until its completion.\footnote{While we restrict to binding and non-preemptive scheduling policies, other common architectures where (a) a server can serve multiple jobs concurrently (processor sharing), (b) a job can be served by multiple servers concurrently, or (c) {job} 
sizes are revealed upon entering the system, are clearly more powerful than the current setting, and are therefore capable of implementing the scheduling policies considered in this paper. As a result, the performance upper bounds developed in this paper also apply to these more powerful variations.} Formally, if a server $j$ has just completed the service of a previous job at time $t$ or is idle, its available actions are: {\bf (a)~Serve a new job}: Server $j$ can choose to fetch a job from any queue in $\neig{j}$ and immediately start service. The server will remain occupied and take no other actions until the processing of the current job is completed, which will take an amount of time that is equal to the size of the job. {\bf (b)~Remain idle}: Server $j$ can choose to remain idle. While in the idling state, it will be allowed to initiate a service (Action $(a)$) at any point in time. 

Given the limited set of actions available to the server, the performance of the system is fully determined by a \emph{scheduling policy}, $\pi$, which specifies for each server $j\in J$,  (a) when to remain idle, and when to serve a new job, and (b) from which queue in $\neig{j}$ to fetch a job when initiating a new service. We only allow policies that are causal, in the sense that the decision at time $t$ depends only on the history of the system (arrivals and service completions) up to $t$. We allow the scheduling policy to be \emph{centralized} (i.e., to have full control over all server actions) based on the knowledge of all \emph{queue lengths} and server states. On the other hand, the policy does \emph{not} observe the actual sizes of the jobs before they are served.

\subsection{Performance Metrics}
\label{sec:ArchiPerMetri}
\bfpara{Characterization of Arrival Rates.} 
{We} will restrict ourselves to arrival rate vectors with average \emph{traffic intensity} {at most} $\rho$, i.e.,
\begin{equation}
\sum_{i=1}^n \lambda_i \leq \rho n,
\end{equation}where $\rho\in(0,1)$ {will be treated throughout the paper as a given absolute constant.} To quantify the \emph{level of variability} or \emph{uncertainty} of a set of arrival rate vectors, $\bLam$, we introduce a \emph{fluctuation parameter}, denoted by $u_n$, {with the property that $\lambda_i<u_n$, for all $i$ and
$\blam\in\bLam$.}

Note that, for a graph with maximum degree $d_n$, the fluctuation parameter should not exceed $d_n$, because otherwise {there could exist some $\lambda\in \bLam$ under which} at least one queue {would} be unstable. Therefore, the best we can hope for is a flexible architecture that can accommodate arrival rate vectors with a $u_n$ that is close to $d_n$. The following condition formally characterizes the range of arrival rate vectors we will be interested in, parameterized by the fluctuation parameter, $u_n$, and traffic intensity, $\rho$. 

\begin{cond} \label{cond:arrdis} {\bf (Rate Condition)} Fix $n\geq 1$ and some $u_n>0$. We say that a (non-negative) arrival rate vector $\blam$ satisfies the rate condition if the following hold: 
\begin{enumerate} 
	\item $\max_{1\leq i \leq \llfl n\rrfl}\lambda_{i} < \lun$. 
	\item $\sum_{i=1}^{n} \lambda_{i} \leq \rho n$.
\end{enumerate}
We denote by $\bLam_n(u_n)$ the set of all arrival rate vectors that satisfy the above conditions.
\end{cond}

\noindent
\bfpara{Capacity Region.} The capacity region for a given architecture is defined as the set of all arrival rate vectors that it can handle. As mentioned in the Introduction, a larger capacity region indicates that the architecture is more robust against uncertainties or changes in the arrival rates. More formally, we have the following definition. 

\begin{definition}[Feasible Demands and Capacity Region]
\label{def:deFlow}
Let ${g} = (I\cup J,E)$ be an $n\times {n'}$ bipartite graph. An arrival rate vector $\blam = \left(\lambda_1,\ldots, \lambda_n\right)$, is said to be feasible if there exists a flow, $
{\mathbf{f}} =\{f_{ij}: (i,j)\in E\}$, such that
\begin{align}
\lambda_i = \sum_{j \in \neig{i}} f_{ij}, &\quad \forall i \in I, \nnb\\
\sum_{i \in \neig{j}} f_{ij} < 1, &\quad \forall j \in J,  \nnb \\
f_{ij} \geq 0,& \quad \forall (i,j)\in E.
\end{align}
In this case, we say that the flow ${\mathbf{f}}$ \emph{satisfies} the demand $\blam$. The \emph{capacity region} of ${g}$, denoted by $\BR({g})$, is defined as the set of all feasible demand vectors of ${g}$. 
\end{definition}

{It is well known that there exists a policy under which the steady-state expected delay is finite if and only if $\blam\in\BR(g_n)$; the strict inequalities in Definition \ref{def:deFlow} are important here.}
For the remainder of the paper, we will use the fluctuation parameter $u_n$ (cf.~Condition \ref{cond:arrdis}) to gauge the size of the capacity region, $\BR(g_n)$, of an architecture. For instance, if $\bLam_n(u_n)\subset \BR(g_n)$, then the architecture $g_n$, {together with a suitable scheduling policy, allows for finite steady-state expected delay, for any arrival rate vetor in $\bLam_n(u_n)$.} 

\noindent
\bfpara{{Vanishing} Delay.} We define the \emph{expected average delay}, $\EX{W|\blam, g, \pi}$ under the arrival rate vector $\blam$, flexible architecture $g$, and scheduling policy $\pi$, {as follows.}  We denote by $W_{i,m}$ the waiting time in queue experienced by the $m$th job arriving to queue $i$, {define
$$\EX{W_i} = \limsup_{m\to \infty}\EX{W_{i,m}},$$}
and  let
\begin{equation}
\EX{W\smid \blam, g,\pi} = \frac{1}{\sum_{i \in I}\lambda_i} \sum_{i\in I}\lambda_i\EX{W_i}.
\end{equation}
  In the sequel, we will often omit the mention of $\pi$, and sometimes of $g$, and write $\EX{W\smid\blam,g}$ or $\EX{W\smid\blam}$, in order to place emphasis on the dependencies that we wish to focus on.\footnote{Note that $\EX{W|\blam}$ captures a \emph{worst-case} expected waiting time across all jobs in the long run, and is always well defined, even under scheduling policies that do not induce a steady-state distribution.}

The delay performance of the system is measured by  the following criteria: (a)~for what ranges $\bLam_n(u_n)$ of arrival rates, $\blam$, 
does delay diminish {to zero} as the system size increases,
i.e., $\sup_{\blam \in\bLam_n(u_n)}\EX{W\smid\blam} \to 0$ as $n \to \infty$, and (b) at what \emph{speed} does the delay diminish, as a function of $n$?

\subsection{Notation}

We will denote by $\N$, $\zp$, and $\rp$, the sets of natural numbers, non-negative integers, and non-negative reals, respectively. The following short-hand notation for asymptotic comparisons will be used {often}, {as an alternative to the usual $\bigo{\cdot}$ notation;} here $f$ and $g$ are positive functions, and $L$ is a certain limiting value of interest, in the set of extended reals, $\R \cup \{-\infty, +\infty\}$: 
\begin{enumerate} 
\item $f(x)\lsim g(x)$ {or $g(x)\gsim f(x)$ for ${\limsup}_{x \to L}{f(x)}/{g(x)}<\infty$;}
\item {$f(x)\ll g(x)$ or $g(x)\gg f(x)$ for ${\limsup}_{x \to L}{f(x)}/{g(x)}=0$;}
\item $f(x)\sim g(x)$ for $\lim_{x \to L}{f(x)}/{g(x)}=1$.
\end{enumerate}

We will minimize the use of floors and ceilings, to avoid the cluttering of notation, and thus assume that all values of interest are appropriately rounded up or down to an integer, whenever doing so does not cause ambiguity or confusion. Whenever {suitable}, we will use upper-case letters for random variables, and lower-case letters for deterministic values. 

\section{Main Results: Capacity Region and Delay of  Flexible Architectures}
\label{sec:mResults}

The statements of our main results are given in this section. Below is a high-level summary of our results; a more complete comparison is given in Table \ref{tab:summary}.

\begin{table}[ht]
\begin{tabular}{|l|l|l|l|}
\hline
\rule{0pt}{2.5ex}  Flexible architectures
&  \parbox{2.8cm}{Rate Conditions}
& Capacity Region
& Delay\\[1ex]
 \hline

\rule{0pt}{6ex}\parbox{4cm}{{\bf Expander} \\(Theorem \ref{thm:expanderArch})} 
& \parbox{2.8cm}{$d_n \gg \ln n$, \\$u_n\lsim d_n $}
& \parbox{3.3cm}{Good for all $\blam$}
& \parbox{4cm}{\vspace{2pt} Good for all $\blam$, with  \\ $\EX{W} \lsim {\ln  n}/{d_n}$, } \\[3.5ex]
 \hline
\rule{0pt}{6ex}\parbox{4cm}{{\bf Modular}\\(Theorems \ref{thm:ModCapacity}, \ref{thm:dimDelayMod})}
& \parbox{2.8 cm}{$d_n\gg 1$, \\$u_n > 1$}
& \parbox{3.3cm}{Bad for {some} \\$\blam$ (even if $u_n\lsim 1$)}
& \parbox{4.2cm}{\vspace{2pt}  Good for uniform $\blam$, with \\{\small $\EX{W} \lsim \exp(-c\cdot d_n)$}} 
\\[3.5ex]
 \hline
\rule{0pt}{7ex}\parbox{4cm}{{\bf Random Modular}\\ {\small (w.h.p.)} \\(Theorems \ref{thm:randModular}, \ref{thm:dimDelayMod})}
& \parbox{2.8 cm}{$d_n\gsim \ln n$, \\$u_n\lsim d_n/\ln n$}
&\parbox{3.3cm}{Good for most $\blam$,\\Bad for some $\blam$}
&\parbox{4cm}{Good for most $\blam$, with \\{\small $\EX{W} {\lsim} \exp(-c\cdot d_n)$},\\ Bad for some $\blam$}  \\[4ex]
 \hline
%
%
%

\end{tabular}

\vspace{20pt}

\caption[justification=justified]
{
\flushleft This table summarizes and compares the flexibility architectures that we study, {in terms of}  of capacity and delay. We say that capacity is ``good'' for $\blam$ if $\blam$ falls within the capacity region of the architecture, and that delay is ``good'' if the expected delay is vanishingly small for large $n$. When describing the size of the set of $\blam$ for which a statement applies, we use the following (progressively weaker) quantifiers: \vspace{5pt} \\ 
{\bf 1.} ``{F}or all'' means that the statement holds for all $\blam\in \bLam_n(u_n)$; \vspace{5pt} \\ 
{\bf 2.} ``{F}or most'' means that the statement holds with high probability when $\blam$ is drawn from an arbitrary distribution over $\bLam_n(u_n)$, independently from any randomization in the construction of the flexibility architecture; \vspace{5pt} \\ 
{\bf 3.} ``{F}or  {some}'' means that the statement is true for a non-empty set of {values of} $\blam$. \vspace{5pt} \\ 
The label ``w.h.p.''  means that all statements in the corresponding row hold with high probability with respect to the randomness in generating the flexibility architecture.
}
\label{tab:summary}
\end{table}

{Our main results focus on {an} Expander architecture, where the interconnection topology is an expander graph with appropriate expansion. We show that, when $d_n \gg \ln n$,  the Expander architecture {has} a capacity region that is within a constant factor of the maximum possible among all graphs with average degree $d_n$, while simultaneously ensuring an asymptotically diminishing queueing delay of order $\ln n/d_n$ for all arrival rate vectors in the capacity region, as $n\to \infty$ (Theorem \ref{thm:expanderArch}). Our analysis {involves} on a new class of virtual-queue-based scheduling policies that rely on dynamically constructed job-to-server assignments on the connectivity graph. }

{Our secondary results concern {a} Modular architecture, which has a simpler construction and scheduling rule compared to the Expander architecture.} The Modular architecture consists of a collection of \emph{separate}  smaller subnetworks, {with complete connectivity between all queues and servers within each subnetwork.} Since the subnetworks are disconnected from {each other}, a Modular architecture does \emph{not} admit a large capacity region: there always exists an \emph{infeasible} arrival rate vector even when the fluctuation  parameter is of constant order (Theorem \ref{thm:ModCapacity}). Nevertheless, we show that with proper randomization in the construction of the subnetworks (Randomized Modular architecture), a simple greedy scheduling policy is able to deliver asymptotically vanishing delay for ``most'' arrival rate vectors with nearly optimal fluctuation parameters, with high probability (Theorem \ref{thm:randModular}). {These findings suggest that, thanks to its simplicity, the Randomized Modular architecture could be a viable alternative to the Expander architecture {if the robustness requirement} is not as stringent and one is content with probabilistic guarantees on system stability.}

\subsection{Preliminaries}
\label{s:prel}

Before proceeding, we provide some information on expander graphs, which will be used in some of our constructions and proofs. 

\begin{definition} An $n\times n'$  bipartite graph $(I\cup J, E)$ is an $(\alpha,\beta)$-expander, if for all $S\subset I$ that satisfy $|S|\leq \alpha n$, we have that $|\neig{S}|\geq \beta|S|$, where $\neig{S} = \bigcup_{i\in S} \neig{{i}}$ is the set of nodes in $J$ that are connected to {some} node in $S$.
\end{definition}

The usefulness of expanders in our context comes from the following lemma, which relates 
the parameters of an expander
to the size of its  capacity region, as measured by the fluctuation parameter, $u_n$. The proof is elementary and is given in Appendix \ref{app:lem:expandRobust}. 

\begin{lemma}[Capacity of Expanders] \label{lem:expandRobust}
Fix $n, n' \in \N$, $\rho\in (0,1)$, $\gamma >\rho$. Suppose that {an $n\times n'$ bipartite graph}, $g_n$, is a {$\left({\gamma}/{\beta_n}, \beta_n\right)$-expander, where $\beta_n\geq u_n$}. Then  $\bLam_n(u_n) \subset \BR(g_n)$. \end{lemma}

The following lemma ensures that such expander graphs exist for the range of parameters that we are interested in. The lemma is a simple consequence of a standard result on the existence of expander graphs, and its proof is given in Appendix \ref{app:cor:expan}. 

\begin{lemma}
\label{cor:expan}
Fix $\rho \in (0,1)$. Suppose that $d_n\to \infty$ as $n \to \infty$. Let $\beta_n = \frac{1}{2} \cdot \frac{\ln (1/\rho)}{1+\ln (1/\rho)} \ d_n$, and $\gamma = \sqrt{\rho}$. There exists $n'>0$, such that for all $n\geq n'$, there exists an $n\times n$ bipartite graph which is a $(\gamma/\beta_n\, ,\beta_n)$-expander with maximum degree $d_n$. 
\end{lemma}

{\noindent
\bfpara{Remark.} It is well known that random graphs with appropriate {average degree} are expanders with high probability (cf.~\cite{hoory2006expander}). For instance, it is not difficult to show that if $d_n\gg \ln n$ and $\beta_n = \frac{1-\gamma}{4}d_n/\ln n$, then an Erd\"os-R\'enyi random bipartite graph with average degree $d_n$ is a $(\gamma/\beta_n, \beta_n)$-expander, with high probability, as $n\to \infty$ (cf.~Lemma 3.12 of \cite{xu2014power}). We note, however, that to {deterministically} construct expanders in a computationally efficient manner can be challenging and is in and of itself an active field of research; the reader is referred to the survey paper \cite{hoory2006expander} and the references therein.
}

\subsection{Expander Architecture}
\label{sec:RandomGraphArchi}

\bfpara{Construction of the Architecture.} {The connectivity graph in the Expander Architecture is an expander graph with maximum degree $d_n$ and appropriate  expansion.} 

\noindent
\bfpara{Scheduling Policy.} 
We employ a  scheduling policy that organizes the arrivals into batches, stores the batches in a virtual queue, and dynamically assigns the jobs in a batch to appropriate servers. 
{Theorem \ref{thm:expanderArch},
which is the main result of this paper,  shows  that under this policy} the Expander architecture achieves an asymptotically \emph{vanishing} delay for \emph{all} arrival rate vectors in the set $ \bLam_n(u_n)$. Of course we assume that $d_n$ is sufficiently large so that the corresponding expander graph exists (Lemma \ref{cor:expan}, with $\rho$ replaced with $\hat{\rho}$). {At a high level, the strong guarantees stem from the excellent connectivity of an expander graph, and similarly of random subsets of an expander graph, a fact which we will exploit to show that jobs arriving to the system during a small time interval can be quickly assigned to connected idle servers with high probability, which then leads to a small delay.} The proof of the theorem, including a detailed description of the scheduling policy, is given in Section \ref{sec:VirtualQueues}.

\begin{theorem}[Capacity and Delay of Expander Architectures]\label{thm:expanderArch}  Let $\hat{\rho}=\frac{1}{1+(1-\rho)/8} $. For every $n\in \N$, define
 $$\beta_n = \frac{1}{2} \cdot \frac{\ln(1/\hat{\rho})}{\ln(1/\hat{\rho})+1} \ d_n, $$ and 
$$\gamma = {\sqrt{\hat{\rho}}}.$$
Suppose that $\ln n \ll d_n\ll n$, and $$u_n\leq \frac{1-\rho}{2} \beta_n.$$ 
Let $g_n$ be a $(\gamma/\beta_n, \beta_n)$-expander with maximum degree $d_n$. The following holds. 
\begin{enumerate}
\item There exists a scheduling policy, $\pi_n$, under which
\begin{equation}
\sup_{\blam_n \in \bLam_n(u_n)} \EX{W | \blam_n, g_n} \leq \frac{c \ln n}{ d_n},
\label{eq:expArchDelay}
\end{equation}
where $c$ is a constant independent of $n$ and $g_n$. 

\item The scheduling policy, $\pi_n$, only depends on $g_n$ and an upper bound on the traffic intensity, $\rho$. It does not require knowledge of the arrival rate vector $\blam_n$. 
\end{enumerate}
\end{theorem}

Note that when $\rho $ is viewed as a constant, the upper bound on $u_n$ in the statement of Theorem \ref{thm:expanderArch} is just a constant multiple of $d_n$. Since the fluctuation parameter, $u_n$, should be no more than $d_n$ for stability to be possible, the size of $\bLam_n(u_n)$ in Theorem \ref{thm:expanderArch} is within a constant factor of the \emph{best possible}.  

{\noindent
\bfpara{Remark.} 
Compared to {our earlier results, in a preliminary version of this paper (Theorem 1 in \cite{TX13SIG}),} Theorem \ref{thm:expanderArch} is stronger in two major aspects: $(1)$ the guarantee for diminishing delay holds deterministically over all arrival rate vectors in $\bLam_n(u_n)$, as opposed to ``with high probability'' {over} the randomness in the generation of $g_n$, and $(2)$ the fluctuation parameter, $u_n$, is allowed to be of order $d_n$ in Theorem \ref{thm:expanderArch}, while \cite{TX13SIG} required that $u_n\ll \sqrt{d_n/\ln n}$. The flexible architecture considered in \cite{TX13SIG}  {was based on Erd\" os-R\' enyi random graphs. It also  employed a scheduling policy based on virtual queues, as in this paper. However, the policy in the present paper is  simpler to describe and analyze.}}

\subsection{Modular Architectures}
\label{sec:modular}

\begin{figure}
\centering
\includegraphics[scale=.5]{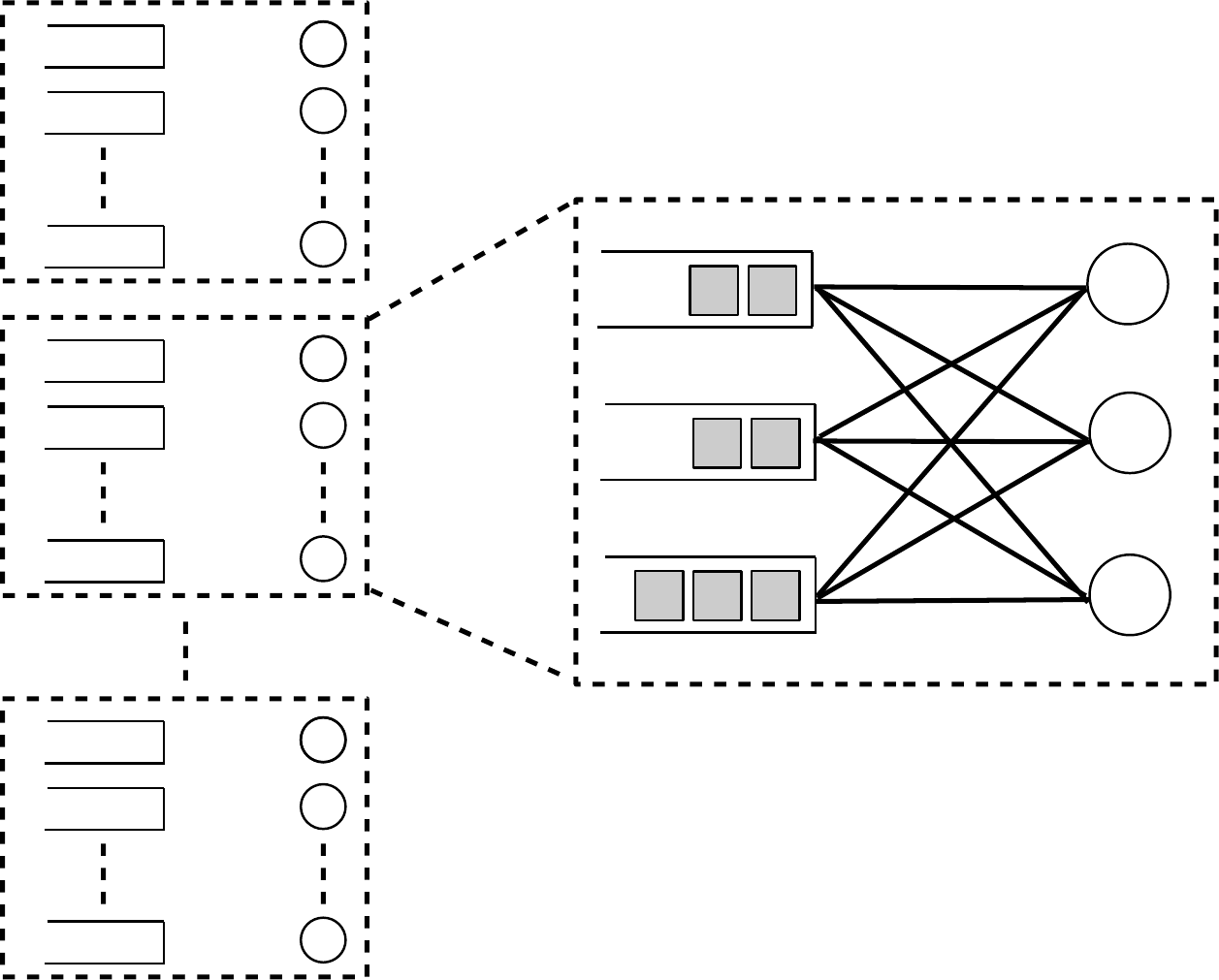}
\caption{A Modular architecture consisting of $n/{\ldn}$ subnetworks, each with $\ldn$ queues and servers. Within each subnetwork, all servers are connected to all queues.}
\label{fig:3}
\vspace{-15pt}
\end{figure}

In a {Modular architecture}, the designer partitions the network into $n/{\ldn}$ \emph{separate} subnetworks. Each subnetwork consists of $\ldn$ queues and servers that are fully connected (Figure \ref{fig:3}), but disconnected from queues and servers in other subnetworks. 

\noindent
\bfpara{Construction of the Architecture.}  {Formally,} the construction is as follows.

\begin{enumerate}
\item We partition  the set $J$ of servers into $n/d_n$ disjoint subsets (``clusters'') $B_1,\ldots,B_{n/d_n}$, all having the same  cardinality $d_n$. For concreteness, we assign the first $d_n$ servers to the first cluster, $B_1$, the next $d_n$ servers to the second cluster, etc.

\item
We form a partition $\sigma_n=(A_1,\ldots,A_{n/d_n})$ of the set $I$ of queues into $n/d_n$ disjoint subsets (``clusters'') $A_k$, all having the same  cardinality $d_n$.

\item To construct the interconnection topology, for $k=1,\ldots,n/d_n$, we connect every 
queue $i\in A_k$ to every server $j\in B_k$. A pair of queue and server clusters with the same index $k$ will be referred to as a subnetwork. 
\end{enumerate}

Note that in a Modular architecture, the degree of each node is equal to the size, $d_n$, of the clusters. 
Note also that different choices of $\sigma_n$ yield isomorphic architectures. When $\sigma_n$ is drawn uniformly at random from the set of all possible partitions of $I$ into subsets of size $n/d_n$,  we call the resulting topology a \emph{Random Modular} architecture. 

\noindent
\bfpara{Scheduling Policy.} We use a simple greedy policy, equivalent to running each subnetwork as an $M/M/d_n$ queue. Whenever a server $j \in B_k$ becomes available, it starts serving a job from any non-empty queue in $A_k$. Similarly, when a job arrives at queue $i\in A_k$, it is immediately assigned to an arbitrary idle server in $B_k$, if such a server exists, and waits in queue $i$, otherwise. 

Our first result points out that a Modular architecture does not have a large capacity region: for any partition $\sigma_n$, there always exists {an} {in}feasible arrival rate vector, even if $u_n$ is small, of order $\mathcal{O}(1)$. The proof is given in Appendix \ref{app:thm:ModCapacity}. Note that this is a negative result that applies no matter what scheduling policy is used.

\begin{theorem}[Capacity Region of Deterministic Modular Architectures] \label{thm:ModCapacity}
Fix $n\geq 1$ {and some $u_n>1$.} 
Let $g_n$ be a Modular architecture {with average degree $d_n\leq \frac{\rho}{2}n$.}  Then,  there exists $\blam_n\in \bLam_n(u_n)$ such that  $\blam_n \notin \BR(g_n)$. 
\end{theorem}

However, if we are willing to 
settle for a weaker result on the capacity region,
 the next theorem states that with the Random Modular architecture,
any given arrival rate vector $\blam_n$ has high probability (with respect to the {random choice of the partition} $\sigma_n$)  of belonging to the capacity region,
if the fluctuation parameter, $u_n$, is of order $\mathcal{O}(d_n / \ln n)$, but no more than that. {Intuitively, this is because the randomization in the connectivity structure makes it unlikely that many large components of $\blam_n$ reside in the same sub-network.} The proof is given in Appendix \ref{app:thm:randModular}.

\begin{theorem}
[Capacity Region of Random Modular Architectures] \label{thm:randModular}
Let $\sigma_n$ be drawn uniformly at random from the set of all partitions, and let $G_n$ be the resulting Random Modular architecture. 
Let $\pb_{G_n}$ be the probability measure that describes the distribution of $G_n$. Fix a constant $c_1>0$, and
suppose that $d_n\geq c_1\ln n$. Then, there exist positive constants $c_2$ and $c_3$, such that:
\begin{enumerate}
\item[(a)] If $u_n\leq c_2d_n/\ln n$, then
\begin{equation}
\lim_{n\to \infty} \inf_{\blam_n \in \bLam_n(u_n)} \pb_{G_n}\left(\blam_n \in \BR(G_n)\right) = 1. 
\label{eq:randModUpper}
\end{equation}
\item[(b)] Conversely,  if $u_n >  c_3d_n/ \ln n$ and $d_n\leq n^{0.3}$, then 
\begin{equation}
\lim_{n\to \infty} \inf_{\blam_n \in \bLam_n(u_n)} \pb_{G_n}\left( \blam_n \in \BR(G_n) \right) = 0,
\label{eq:randModLower}
\end{equation}
\end{enumerate} 
\end{theorem}

We can use {Theorem \ref{thm:randModular}} to obtain a statement about ``most'' arrival rate vectors in $\bLam_n(u_n)$, as follows. Suppose that $\blam_n$ is drawn from an arbitrary distribution $\mu_n$ over $\bLam_n(u_n)$, independently from the randomness in $G_n$. Let $\pb_{G_n}\times \mu_n$  be the product measure that describes the joint distribution of $G_n$ and $\blam_n$. Using Fubini's theorem, Eq.\ \eqref{eq:randModUpper} implies that
\begin{equation}
\lim_{n\to \infty} (\pb_{G_n}\times \mu_n)\big(\blam_n \in \BR(G_n)\big) = 1. 
\label{eq:randModUpper2}
\end{equation}
{A further} application of Fubini's theorem {and an elementary argument}\footnote{{We are using here the following elementary Lemma.
Let $A$ be an event with $\pb(A)\geq 1-\epsilon$,  and let $X$ be a random variable. Then, there exists a set $B$ with $\pb(B)\geq 1-\sqrt{\epsilon}$ such that 
$\pb(A\mid X)\geq 1-\sqrt{\epsilon}$, whenever $X\in B$.
The lemma is applied by letting $A$ be the event $\{\blam_n \in \BR(G_n)\}$ and letting $X=G_n$.}}
implies that there exists a sequence $\delta'_n$ that converges to zero, such that the event 
\begin{equation}
\mu_n\big(\blam_n \in \BR(G_n) \mid G_n\big)\geq 1-\delta'_n
\label{eq:randModUpper3}
\end{equation}
has ``high probability,'' with respect to the measure $\pb_{G_n}$. That is, there is high probability that the Random Modular architecture includes ``most'' arrival vectors $\blam_n$ in $\bLam_n(u_n)$.


We now turn to delay. The {next} theorem states that in a Modular architecture, delay is vanishingly small for all arrival rate vectors in the capacity region that are not too close to its outer boundary. The proof is given in Appendix \ref{app:thm:dimDelayMod}. 

We need some notation. For any set $S$ and scalar $\gamma$, we let $\gamma S=\{\gamma x: x\in S\}$.

\begin{theorem}[Delay of Modular Architectures] \label{thm:dimDelayMod}
Fix some $\gamma\in (0,1)$, and {consider   a Modular architecture $g_n$ for each $n$.}  There exists a constant $c>0$, independent of $n$ {and the sequence $\{g_n\}$}, so that 
\begin{equation}
\EX{W \, |\, \blam_n} \lsim \exp(-c \cdot d_n),
\end{equation}
for every $\blam_n\in \gamma \BR(g_n)$.
\end{theorem}

\subsubsection{Expanded Modular Architectures}

{There is a further} variant of the Modular architecture {that we call} the \emph{Expanded Modular architecture}, which combines the features of a Modular architecture and an expander graph via a graph product. By construction, it 
uses part of the system flexibility to achieve a large capacity region and part to achieve low delay.
As a result, the Expanded Modular architecture admits a smaller capacity region compared to that of {an} Expander architecture. 
{Another drawback is that the available performance guarantees involve policies that require the knowledge of the arrival rates $\lambda_i$.}
On the positive side, it guarantees {an asymptotically vanishing} delay for \emph{all} arrival rates, uniformly across the capacity region, and can be operated by a scheduling policy that is arguably simpler than in the Expander architecture. The construction and {a} scheduling {policy} for the Expanded Modular architecture {is given} in Appendix \ref{app:expanModular}, along with {a statement of its} performance guarantees (Theorem \ref{thm:expandModular}). {The technical details can be found in \cite{xu2014power}.} 

\input{VirtualQueueSection_V15_Rev2-Final}

\section{Summary and Future Research }
\label{sec:discFut}

The main message of this paper is that {the two objectives of} a large capacity region and {an asymptotically vanishing} delay can be {simultaneously} achieved {even if} the level of processing flexibility of each server is small compared to the system size. Our main results show that, 
{as far as these objectives are concerned,} 
the family of Expander architectures is essentially optimal: it admits a capacity region whose size is within a constant factor of the maximum possible, while ensuring an asymptotically vanishing queueing delay for all arrival rate vectors in the capacity region. 

{
An alternative design}, the Random Modular architecture, {guarantees} small delays for ``many'' arrival rates, by means of a simple greedy scheduling policy. However, for any given Modular architecture, there are always many arrival rate vectors {in $\bLam_n(u_n)$ that result in an unstable system,} even if the maximum arrival rate across the queues is of constant order. {Nevertheless, the simplicity of the Modular architectures can still be  appealing in some practical settings.}

Our result for the Expander architecture leaves open {three} questions: 
\begin{enumerate}
\item Is it possible to lower the requirement on the average degree from $d_n\gg \ln n$  to $d_n \gg 1$? 
\item Without sacrificing the size of the capacity region, is it possible to achieve a queueing delay which approaches zero exponentially fast  as a function of $d_n$? The delay scaling in Theorem \ref{thm:expanderArch} is $\mathcal{O}(\ln n / d_n)$. 
\item {Is it possible to obtain delay and stability guarantees under simpler policies, such as the greedy heuristic mentioned in Section \ref{sec:simu}? {The techniques developed in \cite{VADW12} for analyzing first-come-first-serve scheduling rules in a multi-class queueing network similar to ours could be a useful starting point.} }
\end{enumerate}

\del{DELETE.
\begin{conjecture}[Existence of ``Ideal'' Architectures]
\label{conj:perfect}
Suppose that $d_n\gg 1$. There exists a constant $h>0$, such that if  
\begin{equation}
\frac{u_n}{d_n} \leq  h,\quad \mbox{for all $n\geq 1$},
\end{equation} 
then there exists a sequence of architectures, $\{g_n\}_{n\geq 1}$, and associated scheduling policies, under which
\begin{equation}
\EX{W\,|\, g_n,\blam} \leq c_1\exp(-c_2\cdot d_n),\quad \mbox{for all $n\geq 1$ and $\blam\in \bLam_n(u_n)$},
\label{eq:conjeDelayScal}
\end{equation}
where $c_1$ and $c_2$ are positive constants independent of $n$ or $\blam$. 
\end{conjecture}
}


{Finally,} the scaling regime considered in this paper assumes that the traffic intensity is fixed as $n$ increases, which fails to capture {system} performance in the heavy-traffic regime ($\rho\approx1$). {It} would be interesting to consider a scaling regime in which $\rho$ and $n$ scale simultaneously (e.g., as in the celebrated Halfin-Whitt regime \cite{HW81}), but it is unclear at this stage 
what the most appropriate formulations and analytical techniques are.

\bibliographystyle{abbrv}
\bibliography{mainbibv2}

\newpage


%
%
%
%
%

\begin{APPENDICES}

\section{Proofs}

\subsection{Proof of Lemma \ref{lem:expandRobust}}\label{app:lem:expandRobust}

\bpf Fix $\blam = (\lambda_1,\ldots, \lambda_n) \in \bLam_n(u_n)$, and let $g_n$ be a $(\gamma/\beta_n, \beta_n)$-expander, where $\gamma > \rho$ and $\beta_n\geq u_n$. By the max-flow-min-cut theorem, and the fact that all servers have unit capacity,  it suffices to show that \begin{equation}
\sum_{i \in S} \lambda_i < |\neig{S}|, \quad \forall S \subset I. 
\end{equation}
We consider two cases, depending on the size of $S$. \begin{enumerate}
\item Suppose that $|S|< \gamma n/{\beta_n}$. By the expansion property of $g_n$, we have that 
\begin{equation}
\neig{S} {\geq} {\beta_n|S| \geq u_n|S| }  {>} \sum_{i \in S} \lambda_i,
\label{eq:expRob1}
\end{equation}
where the second inequality follows from the fact that $\beta_n\geq u_n$, and the last inequality from $\lambda_i {<} u_n$ for all $i\in I$. 
\item Suppose that $|S| \geq {\gamma}n/{\beta_n}$. 
{By removing, if necessary, some of the nodes in $S$, we obtain a set $S'\subset S$ of size exactly $\gamma n/\beta_n$, and}
\begin{equation}
\neig{S} {\geq \neig{S'}} \stackrel{(a)}{\geq} \gamma n > \rho n \stackrel{(b)}{\geq} \sum_{i \in S}\lambda_i,
\end{equation}
where step $(a)$ follows from the {expansion property}, and step $(b)$ from the assumption that $\sum_{i\in I}\lambda_i \leq \rho n$. 
\end{enumerate}
This completes the proof. 
\qed

\subsection{Proof of Lemma \ref{cor:expan}}
\label{app:cor:expan}
\bpf
Lemma \ref{cor:expan} is a consequence of the following standard result ({cf.} \cite{Asr98}), where we let  $d= d_n$, $\beta=\beta_n$, and $\alpha= \gamma/\beta_n = \sqrt{\rho}/\beta_n$, and observe that $\log_2\beta_n \ll \beta_n$ as $n\to \infty$. 
\begin{lemma} Fix $n\geq 1$, $\beta \geq  1$ and $\alpha \beta <1$. 
{If}
\begin{equation}
d\geq \frac{1+\log_2\beta +(\beta+1)\log_2 e}{-\log_2(\alpha \beta)}+\beta+1,
\end{equation}
{then}
there exists an $(\alpha, \beta )$-expander with maximum degree $d$. 
\end{lemma}
\qed

\subsection{Proof of Theorem \ref{thm:ModCapacity}}
\label{app:thm:ModCapacity}
\bpf
Since the arrival rate vector $\blam_n$ whose existence we want to show  can depend on the architecture, we assume, without loss of generality, that servers and queues are clustered in the same manner: server $i$ and queue $i$ belong to the same cluster. Since all servers have capacity 1, and each cluster has exactly $d_n$ servers, it suffices to show that there exists $\blam = (\lambda_1,\ldots, \lambda_n) \in \bLam_n(u_n)$, such that the total arrival rate to the first queue cluster exceeds $d_n$, i.e.,
\begin{equation}
\sum_{i=1}^{d_n}\lambda_i >d_n.
\label{eq:badlam}
\end{equation}
To this end, consider the vector $\blam$ where $\lambda_i = \min\{2, (1+u_n)/2\} $ for all $i\in \{1,\ldots, d_n\}$, and $\lambda_i = 0$ for $i\geq d_n+1$. Because of the assumption $u_n>1$ in the statement of the theorem, we have that
\begin{equation}
 \max_{1\leq i\leq n}\lambda_i =\min\{2,(1+u_n)/2\}  \leq \frac{1+u_n}{2}<u_n,
 \label{eq:badlam1}
 \end{equation} 
 and 
\begin{equation}
\sum_{i=1}^n  \lambda_i = d_n\min\{2,(1+u_n)/2\} \leq 2d_n \leq  2\cdot\frac{\rho}{2} n = \rho n, \label{eq:badlam2}
\end{equation}
where the last inequality in Eq.\ \eqref{eq:badlam2} follows from the assumption that $d_n\leq \frac{\rho}{2}n$. Eqs.~\eqref{eq:badlam1} and \eqref{eq:badlam2} together ensure that $\blam \in \bLam_n(u_n)$ (cf.\ Condition 1). Since we have assumed that $u_n>1$, we have $\lambda_i>1$, for $i=1,\ldots,d_n$, and therefore 
Eq.~\eqref{eq:badlam} holds for this $\blam$. We thus have that $\blam \notin \BR(g_n)$, which proves our claim. 
\qed

\subsection{Proof of Theorem \ref{thm:randModular}}
\label{app:thm:randModular}

\bpf 
{\bf Part (a); Eq.~\eqref{eq:randModUpper}}. We will use the following classical result due to Hoeffding, adapted from Theorem $3$ in \cite{Hoeff63}.

\begin{lemma}
\label{lem:hoeff}
Fix integers $m$ and $n$, where $0<m<n$. Let $X_1,X_2,\ldots, X_m$ be random variables drawn uniformly  from a finite set $C = \left\{c_1,\ldots, c_n\right\}$, without replacement. Suppose that $0\leq c_i \leq b$ for all $i$, and let $\sigma^2=\var{X_1}$. {Let} $\overline{X} = \frac{1}{m}\sum_{i=1}^m X_i$. {Then,}
\begin{equation}
	\pb\left(\overline{X} \geq \EX{\overline{X}} + t \right) \leq \exp\left(-\frac{mt}{b}\left[\left(1+\frac{\sigma^2}{bt}\right)\ln\left(1+\frac{bt}{\sigma^2}\right)-1\right]\right),
\end{equation}
for all $ t\in (0,b)$.
\end{lemma}

We fix some $\blam_n \in \bLam_n(u_n)$.  If $u_n<1$, then {$\blam_n \in \bLam_n(1)$. It therefore suffices to prove the result for the case where $u_n\geq 1$ and we will henceforth assume that this is the case.} 
Recall that $A_k \subset I$ is the set of $d_n$ queues in the $k$th queue cluster generated by the partition $\sigma_n=(A_1,\ldots,A_{n/d_n})$. 
We consider some $\epsilon\in(0,1/\rho)$, 
and define the event $E_k$ as 
\begin{equation}
	E_k = \left\{\sum_{i\in A_k} \lambda_i > (1+\epsilon)\rho d_n\right\}.
	\label{eq:EkDef}
\end{equation}
Since $\sigma_n$ is drawn uniformly at random from all possible partitions, it is not difficult to see that  $\sum_{i\in A_k} \lambda_i$ has the same distribution as $ \sum_{i=1}^{d_n}X_i$, where $X_1,X_2,\ldots,X_{d_n}$ are $d_n$ {random variables} drawn uniformly at random, without replacement, from the set $\{\lambda_1,\lambda_2,\ldots,\lambda_n\}$. Note that 
$\epsilon\rho<1\leq u_n$, so that $ \epsilon\rho\in (0,u_n)$. We can therefore apply Lemma \ref{lem:hoeff}, with $m=d_n$, $b=u_n$, and $t=\epsilon\rho$, to obtain
\begin{align}
\pb\left(E_1\right) =& \pb\left(\sum_{i=1}^{d_n}X_i > (1+\epsilon)\rho d_n\right) \nnb\\
\stackrel{(a)}{\leq}& \pb\left(\frac{1}{d_n}\sum_{i=1}^{d_n}X_i > \EX{\frac{1}{d_n}\sum_{i=1}^{d_n}X_i}+ \epsilon\rho \right) \nnb\\
\leq & \exp\left(-\frac{\epsilon \rho d_n}{u_n}\left[\left(1+\frac{\var{X_1}}{\epsilon \rho u_n}\right)\ln\left(1+\frac{\epsilon \rho u_n}{\var{X_1}}\right)-1\right]\right),
\label{eq:hoefinter1}
\end{align}
where the probability is {taken} with respect to the randomness in $G$, and where in step $(a)$ we used the fact that 
\begin{equation}
	\EX{\sum_{i=1}^{d_n}X_i} = \sum_{i=1}^{d_n} \EX{X_i} = d_n\EX{X_1} = d_n\left(\frac{1}{n}\sum_{i=1}^n \lambda_i\right) \leq \rho d_n.
\end{equation}

We now develop an upper bound on $\var{X_1}$.  Since $X_1$ takes values in $[0,u_n]$,
we have $X_1^2 \leq u_n X_1$ and, therefore,
\begin{equation}
	\var{X_1} \leq \E(X_1^2) \leq u_n \E(X_1) \leq \rho u_n. 
	\label{eq:varX1}
\end{equation}
{Observe that for all $a,x>0$,
\begin{equation}
\frac{d}{d x}(1+x/a)\ln(1+a/x)=-\frac{1}{x}+\frac{1}{a}\ln(1+a/x) < -\frac{1}{x}+\frac{1}{a}\cdot \frac{a}{x}=0.
\end{equation}
Therefore, with the substitutions $a=\epsilon\rho u_n$ and $x=\var{X_1}$, we have that the right-hand-side of \eqref{eq:hoefinter1} is increasing in $\var{X_1}$.} Combining Eqs.~\eqref{eq:hoefinter1} and \eqref{eq:varX1},  we obtain
\begin{align}
\pb\left(E_1\right) \leq &  \exp\left(-\frac{\epsilon \rho d_n}{u_n}\left[\left(1+\frac{1}{\epsilon}\right)\ln\left(1+{\epsilon}\right)-1\right]\right). \nnb \end{align}
{Note that
\begin{equation}
\frac{d}{dx}\left(1+\frac{1}{x}\right)\ln(1+x) = \frac{1}{x^2}(x-\ln(1+x)) \sk{a}{\to} \frac{1}{2}, \quad \mbox{as $x\downarrow 0$}, 
\end{equation}
where step $(a)$ follows from applying {l}'H$\hat{\mbox{o}}$pital's rule.} We thus have that $\left[\left(1+\frac{1}{\epsilon}\right)\ln(1+\epsilon)-1 \right] \sim \frac{1}{2}\epsilon {\geq} \frac{1}{3}\epsilon,$ as $\epsilon\downarrow 0$, it
follows that there exists $\theta>0$ such that for all $\epsilon\in(0,\theta)$, 
\begin{align}
\pb\left(E_1\right) \leq & \exp\left(-\frac{\rho}{3}\cdot\frac{\epsilon^2 d_n}{u_n}\right).
\label{eq:E1bound}
\end{align}

Let $\epsilon=\frac{1}{2}\min\{\frac{1}{\rho}-1, \theta\}$; in particular, our earlier assumption that $\epsilon\rho<1$ is satisfied. Suppose that $u_n\leq \frac{\rho \epsilon^2}{6} d_n\ln^{-1}n$. Combining Eq.~\eqref{eq:E1bound} with the union bound, we have that
\begin{align}
\pb_{G_n}\left(\blam_n \notin \BR(G_n)\right)	\leq & \pb\left(\bigcup_{k=1}^{n/d_n} E_k\right) \nnb\\
\leq & \sum_{k=1}^{n/d_n} \pb\left(E_k\right) \nnb\\
\leq & \frac{ n}{d_n}\exp\left(-\frac{\rho}{3}\cdot\frac{\epsilon^2 d_n}{u_n}\right) \nnb\\
\stackrel{(a)}{\leq} & \frac{ n}{d_n}\cdot\frac{1}{n^2} \nnb\\
\leq & n^{-1},
\label{eq:lamRGUpper}
\end{align}
where step $(a)$ follows from the assumption that $u_n\leq \frac{\rho \epsilon^2}{6} d_n\ln^{-1}n$. 
It follows that
\begin{align}
\lim_{n\to\infty}\inf_{\blam_n\in\bLam_n(u_n)} \pb_{G_n}\left(\blam_n \in \BR(G_n)\right)	\geq & \lim_{n\to\infty} \Big(1-\frac{1}{n}\Big)=1.
\end{align}
We have therefore proved part (a) of the theorem, with  $c_2={\rho \epsilon^2}/{6}$.

{\bf Part (b); Eq.~\eqref{eq:randModLower}}. 

Let us fix a large enough constant $c_3$, whose  value will be specified later, and let
\begin{equation}
v_n = c_3\frac{d_n}{\ln n}.
\end{equation}
For this part of the proof, we will assume that
$u_n > {v_n}$. Because we are interested in showing a result for the worst case over all $\blam_n \in \bLam_n(u_n)$, we can assume that $u_n\ll n$. 

At this point, we could analyze the model for a worst-case choice of $\blam_n$. However, the analysis turns out to be simpler if we employ the probabilistic method. Denote by $\mu_n$ a probability measure over $\bLam_n(u_n)$. Let $\blam_n$ be a random vector drawn from the distribution $\mu_n$, independent of the randomness in the Random Modular architecture, $G$. (For convenience, we suppress the subscript $n$ and write $G$ instead of $G_n$.)
The following elementary fact captures the essence of the probabilistic method.

\begin{lemma}
\label{lem:MeasureToPoinLam}
Fix n, a measure $\mu_n$ on $\bLam_n(u_n)$, and a constant $a_n$. Suppose that
\begin{equation}
\pb_{\blam_n,  G}\left(\blam_n \notin \BR(G) \right) \geq a_n,  
\label{eq:munanto1}
\end{equation}
where $\pb_{\blam_n,  G}$ stands for the product of the measures $\mu_n$ (for $\blam_n$) and $\pb_{G}$ (for $G$). Then, 
\begin{equation}
\sup_{\tilde{\blam}_n \in \bLam_n(u_n)} \pb_{G}(\tilde{\blam}_n\notin \BR(G))\geq a_n. 
\end{equation}
\end{lemma}

\bpf
We have that
\begin{align}
\sup_{\tilde{\blam}_n \in \bLam_n(u_n)}\pb_{G}(\tilde{\blam}_n \notin \BR(G)) \geq & \int_{\tilde{\blam}_n \in \bLam_n(u_n)} \pb_{G}(\tilde{\blam}_n \notin \BR(G)) \, d\mu_{n}(\tilde{\blam}_n )\nnb\\
= & \pb_{\blam_n,  G}\left(\blam_n \notin \BR(G) \right)\nnb\\
\geq& a_n. 
\end{align}
\qed

We  will now construct sequences, $\{\mu_n: n \in \N\}$, and $\{a_n: n\in \N\}$, with  $\lim_{n \to \infty} a_n=1$, so that Eq.~\eqref{eq:munanto1} holds for all $n$. To simplify notation, in the rest of this proof we will write $\pb$ instead of $\pb_G$ or $\pb_{\blam_n,  G}$, etc. Which particular measure we are dealing with will always be clear from the context.

Fix $n\in \N$. We first construct the distribution $\mu_n$. Let $\blam'=(\lambda_1',\lambda_2',\ldots, \lambda_n' )$ be a random vector with independent components and with 
	\begin{equation}
	\lambda_i' = \left\{
		\begin{array}{l l }
		{v_n}, & \quad  \mbox{w.p. } 		\frac{\rho}{(1+\epsilon){v_n}} ,\\
		0, & \quad \mbox{ otherwise,}
		\end{array}
		\right.
		\label{eq:lamiid}
\end{equation}
for all $i$. Let $H$ be the event defined by
\begin{equation}
H = \left\{\sum_{i=1}^n \lambda'_i \leq \rho n\right\}.
\end{equation}
Let $\blam_n$ be the random vector given by 
\begin{equation}
\blam_n = \iden(H) \blam', 
\end{equation}
where $\mathbf{0}$ is the zero vector {of  dimension $n$,}  and where $\iden(\cdot)$ is the indicator function. That is, $\blam_n$ takes on the value of $\blam'$ if $H$ occurs, and is set to zero, otherwise. It is not difficult to verify that, by construction, we always have $\blam_n\in  \bLam_n(u_n)$. We let $\mu_n$ be the distribution of this random vector $\blam_n$. 

We next show that  
\begin{equation}
\lim_{n \to \infty} \pb\left(\blam_n \notin \BR(G) \right) =1, 
\label{eq:munanto2}
\end{equation}
which, together with Lemma \ref{lem:MeasureToPoinLam} above, will complete the proof of the theorem. Fix some $\epsilon>\frac{1}{\rho}-1$, so that $(1+\epsilon)\rho>1$, and define the event
\begin{align}
E_k = \left\{\sum_{i\in A_k} \lambda'_i > (1+\epsilon)\rho d_n\right\}, \qquad k\in \{1, \ldots, n/d_n\}. 
\end{align}
Note that, if some $E_k$ occurs, then $\blam'$ will not be in $\R(G)$. Therefore,
\begin{equation}
\pb(\blam'\notin \BR(G))\geq \pb\left(\bigcup_{k=1}^{{n}/{d_n}} E_k\right). 
 \end{equation}
Let {$X_1,X_2,\ldots$} be i.i.d.~Bernoulli random variables with 
\begin{equation}
\E(X_1)=\pb(X_1=1) = \frac{\rho }{(1+\epsilon){v_n}}.
\end{equation}
By the definition of $\blam'$ (cf.~Eq.~\eqref{eq:lamiid}), we have that
\begin{align}
\pb(E_1) =& \pb\left(\sum_{i\in A_1} \lambda'_i > (1+\epsilon)\rho d_n\right)\nnb\\
=&\pb\left(\sum_{i=1}^{d_n}X_i > (1+\epsilon)\rho \frac{d_n}{{v_n}}\right) \nnb\\
 = & \pb\left(\frac{1}{d_n}\sum_{i=1}^{d_n}X_i > (1+\epsilon)^2\EX{X_1}\right). 
\end{align}
By Sanov's theorem (cf.~Chapter 12 of \cite{cover2012elements}), we have that
\begin{align}
\pb(E_1) = & \pb\left(\frac{1}{d_n}\sum_{i=1}^{d_n}X_i > (1+\epsilon)^2\EX{X_1}\right) \nnb\\
 \gsim & \frac{1}{d_n^2}\exp\left(- D_B\left( \left. \frac{(1+\epsilon)\rho }{{v_n}} \right\| \frac{\rho }{(1+\epsilon){v_n}}\right)d_n\right),
 \label{eq:PE1iid}
\end{align}
where $D_B(p\| q)$ is the Kullback-Leibler divergence between two Bernoulli distributions with parameters $p$ and $q$, respectively: 
\begin{equation}
D_B(p\|q) = p\ln\frac{p}{q}+(1-p)\ln\frac{1-p}{1-q}. 
\end{equation}

Let us fix some $r \in (0,1)$. Using the fact that $\ln(1+y)\sim y$ as $y\to 0$, we have that
\begin{equation}
D_B\left(\left. x \right \|r x \right) \sim x\left[\ln \frac{1}{r}+(1-r)\right], \quad \mbox{as $x\to 0$}. 
\label{eq:Dxrx}
\end{equation}
Recall that $d_n\geq c_1\ln n$ and ${v_n}\geq {/\ln n}$. By Eq.~\eqref{eq:Dxrx}, with $x =  {(1+\epsilon)\rho}/{{v_n}}$, $r = {1}/{(1+\epsilon)^2}$, and for {the} given $c_1$, we can set $c_3$ to be sufficiently large so that 
\begin{align}
D_B\left( \left. \frac{(1+\epsilon)\rho }{{v_n}} \right\| \frac{\rho }{ (1+\epsilon) {v_n}}\right) \leq & 2\frac{(1+\epsilon)\rho}{{v_n}} \cdot \left[\ln(1+\epsilon)^2+\left(1-\frac{1}{(1+\epsilon)^2}\right)\right] \nnb\\
= & \frac{2h}{{v_n}}, 
\label{eq:D_B2}
\end{align}
for all sufficiently large $n$, where $h = (1+\epsilon)\rho \left[\ln(1+\epsilon)^2+\left(1-\frac{1}{(1+\epsilon)^2}\right)\right] >0$. Combining Eqs.~\eqref{eq:PE1iid} and \eqref{eq:D_B2},  we have that
\begin{align}
\pb(E_1) \gsim \frac{1}{d_n^2}\exp\left(-2h\frac{d_n}{{v_n}}\right) \stackrel{(a)}{\gsim} \frac{1}{d_n^2} n^{-2h/c_3},
\label{eq:E1iidFinal}
\end{align}
where step $(a)$ follows from the assumption that ${v_n} \geq c_3d_n/\ln n$. 
{Equation \eqref{eq:E1iidFinal} can be rewritten in the form
\begin{align}
\pb(E_1) \geq\frac{c}{d_n^2} n^{-2h/c_3},
\label{eq:E1iidFinal2}
\end{align}
where $c$ is a positive constant, and where the inequality is valid for large enough $n$.}

Fix $c_3 = 40h$, and recall that $\epsilon>\frac{1}{\rho}-1$. We have that
\begin{align} 
\pb(\blam' \notin \BR(G))\geq & \pb\left(\bigcup_{k=1}^{{n}/{d_n}} E_k\right) \nnb\\
\stackrel{(a)}{=} & 1-\prod_{k=1}^{n/d_n}\left(1-\pb(E_k)\right)\nnb\\
= & 1- \left(1- \pb(E_1)\right)^{n/d_n}\nnb\\
\stackrel{(b)}{\geq}& 1-\left(1-{c}d_n^{-3}n^{1-2h/c_3}{d_n/n}\right)^{n/d_n}\nnb\\
\stackrel{(c)}{\geq}& 1-\left(1-{c}n^{0.05}{d_n/n}\right)^{n/d_n} \nnb\\
\to& 1, \quad \mbox{as $n\to \infty$,}
\label{eq:iidisbad}
\end{align}
where step $(a)$ is based on the independence among the events {$E_k$}, which is in turn based on the independence among the $\lambda'_i$s; step $(b)$ follows from Eq.~\eqref{eq:E1iidFinal2} and some rearrangement;  step $(c)$ follows from the assumption in the statement of the theorem that $d_n\leq n^{0.3}$,  and {our choice of} $c_3=40h$. 

We next show that the event $H$ occurs with high probability when $n$ is large. Let, as before, the $X_i$s be i.i.d.~Bernoulli random variables with $\E(X_1)=\frac{\rho}{{v_n}(1+\epsilon)}$. Then, 
\begin{align}
\pb(H) =& \pb\left(\sum_{i=1}^n \lambda'_i \leq \rho n\
\right) \nnb\\
= &  \pb\left(\sum_{i=1}^n X_i \leq \rho n/{v_n}\
\right) \nnb\\
= & \pb\left(\frac{1}{n}\sum_{i=1}^{n}X_i \leq (1+\epsilon)\EX{X_1}\right) \to 1, \quad \mbox{as $n\to \infty$}, 
\label{eq:HisBig}
\end{align}
by the weak law of large numbers. 

We are now ready to prove Eq.~\eqref{eq:munanto2}. We have that 
\begin{align}
\pb_{\blam_n, G} \Big( \blam_n \notin \BR(G) \Big) =& \pb_{\blam', G} \Big( \iden(H)\blam' \notin \BR(G) \Big) \nnb\\
=& \pb_{\blam', G} \Big( H \cap \left\{\blam' \notin \BR(G)\right\} \Big) \nnb\\
\geq & \pb(H)+\pb\big( \blam' \notin \BR(G) \big)-1 \nnb\\
\to & 1, \quad \mbox{as $n\to \infty$,}
\label{eq:measureBad}
\end{align}
where the last step follows from Eqs.~\eqref{eq:iidisbad} and \eqref{eq:HisBig}. By Lemma \ref{lem:MeasureToPoinLam}, Eq.~\eqref{eq:measureBad} implies that $\lim_{n\to \infty} \sup_{\blam_n \in \bLam_n(u_n)} \pb_{G_n}\left(\blam_n \notin \BR(G)  \right)=1$, which is in turn equivalent to $\lim_{n\to \infty} \inf_{\blam_n \in \bLam_n(u_n)} \pb_{G_n}\left(\blam_n \in \BR(G)  \right)=0$. This proves Eq.~\eqref{eq:randModLower}. \qed

\subsection{Proof of Theorem 
\ref{thm:dimDelayMod}}
\label{app:thm:dimDelayMod}

\bpf
Denote by $Q_i(t)$ the number of jobs in queue $i$ at time $t$, and by $Q_k(t)$ the total number of jobs in queue cluster $k$, i.e.,  
\begin{equation}
Q_k(t) = \sum_{i\in A_k}Q_i(t). 
\end{equation}
We note that $Q_k(\cdot)$ is the number of jobs in an $M/M/c$ queue, with $c=d_n$ and arrival rate $\eta_k = \sum_{i\in A_k} \lambda_i$. Also note that since $\blam_n \in \gamma \BR(g_n)$, we have that $\eta_k\leq \gamma d_n$. Using the formula for the expected waiting time in queue for an $M/M/c$ queue (cf.~Section 2.3 of \cite{gross2008fundamentals}), one can show that the average waiting time across jobs arriving to cluster $k$, $W_k$, satisfies
\begin{equation}
\EX{W_k | \blam} = \frac{1}{\sum_{i\in A_k}\lambda_i }\sum_{i\in A_k}\lambda_i \EX{W_i}= \frac{C(d_n, \eta_k)}{d_n-\eta_k}\leq  \frac{C(d_n, \gamma d_n)}{(1-\gamma)d_n} {\lsim } \exp(-b \cdot d_n), 
\label{eq:W_qisgoodforMMK}
\end{equation}
where $C(c,r)$ is given by 
{$$C(c,r) = \frac{r^c}{c!}\cdot \frac{1}{c(1-r/c)^2} \left(\frac{r^c}{c!}\cdot \frac{1}{1-r/c}+\sum_{i=0}^{c-1}\frac{r^i}{i!}\right)^{-1}. 
$$} The last inequality in  Eq.~\eqref{eq:W_qisgoodforMMK} follows from the fact that for {any given} $\gamma\in (0,1)$, there exists $b>0$, so that $C(x,\gamma x)\lsim  \exp(-b \cdot x)$ as $x\to \infty$, as can be checked through elementary algebraic manipulations.
\qed

\subsection{Lower Bound on the Total Arrival Rate}
\label{app:LowBonArrival}

We show in this section that the assumption that $\rho \in (1/2 \, , 1)$ and $\sum_{i=1}^n\lambda_i \geq (1-\rho)n$ (cf.~Eq.~\eqref{eq:sumLamLowerBound} in Assumption \ref{ass:lamLowerBound}) can be made without loss of generality. Fix {the} traffic intensity $\rho\in (0,1)$, and suppose that $\blam \in  \bLam_n(u_n)$. Define 
\begin{equation}
\rho' = \rho+\frac{1}{2}(1-\rho)=\frac{1+\rho}{2}. 
\end{equation}
Note that $ 1/2 <\rho'< 1$, and $1-\rho' = (1-\rho)/2.$ Consider a modified vector $\blam'$, where $\lambda'_i =  (1-\rho')+\lambda_i$, for all $i\in \{1,\ldots, n\}$. By construction, we have that
\begin{align}
\sum_{i=1}^n \lambda'_i \geq&  (1-\rho')n, \label{eq:lamPrep0} \\
\sum_{i=1}^n \lambda'_i \leq& (1-\rho')n+\sum_{i=1}^n\lambda_i \leq (1-\rho')n+\rho n = \rho' n, \label{eq:lamPrep1} \\
\max_{1\leq i \leq n}\lambda'_i\leq &  \max_{1\leq i \leq n}\lambda_i +(1-\rho') {<} u_n+(1-\rho'). 
\label{eq:lamPrep2}
\end{align}
The above definition of $\blam'$ amounts to the following: we  feed  each queue with an additional  independent Poisson stream of artificial (dummy) jobs of rate $1-\rho'$. By Eqs.~\eqref{eq:lamPrep1} and \eqref{eq:lamPrep2}, the resulting arrival rate vector, $\blam'$, will belong to the set $\bLam_n(u_n +1-\rho')$. Also, by Eq.~\eqref{eq:lamPrep0}, it will satisfy the lower bound \eqref{eq:sumLamLowerBound} on the total arrival rate, albeit with a modified traffic intensity of $\rho'\in (1/2 \, , \,1)$. Therefore, our assumption can always be satisfied by the insertion of dummy jobs. Note that the increment of $1-\rho'$ to the value of $u_n$ is insignificant in our regime of interest, where $u_n \gg 1$, and the insertion of dummy jobs only requires knowledge of the original traffic intensity, $\rho$.

\subsection{Proof of Lemma \ref{lem:nice} }
\label{app:lem:nice}

\bpf Note that because there are $\rho b_n$ jobs in a batch, the size of $\Gamma$ is at most $\rho b_n$, which is in turn less than $m_n$. This guarantees that the cardinality of $\hat \Gamma$ can be taken to be $m_n$. It {therefore} suffices to show that 
\begin{equation}
\pb\lt(\max_{1\leq i \leq n} A_i \geq  \hat{u}_n \rt) \leq 1/n^3. 
\end{equation}

There is a total of $\rho b_n$ arriving jobs in a single batch, and for each arriving job 
\begin{equation}
\pb\left(\mbox{{the} job arrives to queue $i$}\right) = \frac{\lambda_i}{\sum_{{i=1}}^n\lambda_i} \sk{a}{\leq} \frac{\lambda_i}{(1-\rho) n } \leq \frac{u_n}{(1-\rho) n } \sk{b}{\leq}  \frac{1}{2n}\beta_n \leq \frac{1}{2n\hat{\rho}}\beta_n,
\label{eq:PbjobtoQueuei}
\end{equation}
for all $i$, where steps $(a)$ and $(b)$ follow from the assumptions that $\sum_{{i=1}}^n\lambda_i \geq (1-\rho) n$ (Eq.~\eqref{eq:sumLamLowerBound} in Assumption \ref{ass:lamLowerBound}) and that $u_n\leq \frac{1-\rho}{2} \beta_n$ {(in the statement of Theorem \ref{thm:expanderArch})}, respectively. From Eq.~\eqref{eq:PbjobtoQueuei}, $A_i$ is stochastically dominated by a binomial random variable $\tilde{A}\stackrel{d}{=}\mbox{Bino}(\rho b_n\, , \frac{1}{2n\hat{\rho}}\beta_n)$, with
\begin{align}
  \EX{\tilde{A}} = {\rho} b_n\frac{1}{2n\hat{\rho}}\beta_n =\frac{1}{2}\lt(\beta_n \frac{\rho b_n/\hat{\rho}}{n}\rt)=\frac{1}{2}\lt(\beta_n \frac{m_n}{n}\rt)= \frac{1}{2} \, \hat u_n. 
  \label{eq:EtidA}
    \end{align}
    
Based on this expression of $\EX{\tilde A}$, we will now use an exponential tail bound to bound the probability of the event $\{ \max_{1\leq i \leq n} A_i \geq  \hat{u}_n \}$. 
Recall that $b_n = {\frac{320}{(1-\rho)^2}\cdot \frac{n\ln n}{\beta_n}}$. Using the union bound, we have that 
\begin{eqnarray}
\pb\lt(\max_{1\leq i \leq n} A_i \geq  \hat{u}_n\rt) &= & \pb(A_i \geq \hat u_n, \mbox{ for some }i) \nnb\\
&{\leq}&  {n\pb(A_1 \geq \hat u_n)}\\
&{\leq} & {n} \pb\left(\tilde{A} \geq \hat u_n \right) \nln
&\sk{a}{=} & n\pb\left(\tilde{A} \geq 2 \E(\tilde{A})\right) \nln
&\sk{b}{{\leq}} & n \exp\left(-\frac{1}{3}\E(\tilde A)\right) \nln
&= & n \exp\left(- \frac{\rho}{6\hat{\rho}} \cdot\frac{b_n\beta_n}{n}\right) \nln
&\leq  & n \exp\left(- \frac{\rho}{6} \cdot\frac{b_n\beta_n}{n}\right) \nln
&= & {n} \exp\left(-\frac{\rho }{6}\cdot \frac{{320}}{(1-\rho)^2}\cdot \frac{n\ln n}{\beta_n}\cdot \frac{\beta_n}{n}\right) \nln
&\stackrel{(c)}{\leq} &  {n} \exp\left(-\frac{{160}}{6}\ln n\right) \nln
&{\leq} & n^{-3}.
\label{eq:pbBound}
\end{eqnarray}
Step {$(a)$} follows from Eq.~\eqref{eq:EtidA}. Step $(b)$ {follows} from {the following} multiplicative form of the Chernoff bound (cf.~Chapter 4 of  \cite{mitzenmacher2005probability}), {with $\delta=1$}: $\pb({\tilde A}\geq (1+\delta)\mu) \leq \exp(-\frac{\delta^2}{2+\delta}\mu)$,  where {$\tilde A$} is a binomial random variable with $\E({\tilde A})=\mu$. Step $(c)$ follows from the 
assumption  $\rho\in(1/2,1)$ (cf.\ Assumption \ref{ass:lamLowerBound}), and hence 
\begin{equation}
\frac{\rho }{(1-\rho)^2} \geq \rho \geq 1/2. 
\end{equation}
This completes the proof of Lemma \ref{lem:nice}.  \qed

\subsection{Proof of Lemma \ref{lem:subGraphExpansion}}
\label{app:lem:subGraphExpansion}

\bpf  For a set $S\subset \hat{\Gamma}$, denote by $\mathcal{N}^*(S)$ the set of neighbors of $S$ in $\hat{G}$, i.e., $\mathcal{N}^*(S) = \neig{S}\cap \Delta$. To  {prove} Lemma \ref{lem:subGraphExpansion}, we will leverage the fact that the underlying connectivity graph, $g_n$, is an expander graph with appropriate expansion. As a result, most subsets $S\subset \hat{\Gamma}$ have a large {set} of neighbors, $\neig{S}$, in $g_n$. Because each server in $\neig{S}$ belongs {to} $\mathcal{N}^*(S)$ independently, as a consequence of our scheduling policy, we will then use a concentration inequality to show that, with high probability, the sizes of {the sets} $\mathcal{N}^*(S)$ remain sufficiently large. {Using the union bound over the relevant sets $S$, we will finally conclude that
$\hat{G}$ 
has the desired expansion property, with high probability.}

By the definition of {a} $\lt(\gamma /{\hat{u}_n} \, , \hat{u}_n\rt)$-expander, we {are only} interested in {the expansion of} subsets of $\hat{\Gamma} $ {with size less than or equal to} $|\hat{\Gamma}|{\gamma}/\hat{u}_n$. We first verify below that {the size of such subsets $S$ is sufficiently small to be able to exploit the expansion property of $g_n$ and to infer that $\mathcal{N}^*(S)$ is large.}
We have
\begin{align}
\frac{n \gamma / \beta_n}{|\hat{\Gamma}|{\gamma}/\hat{u}_n} =  \frac{n}{|\hat{\Gamma}|}\cdot \frac{\hat{u}_n}{\beta_n} =    \frac{n}{m_n}\cdot \frac{\beta_n\frac{m_n}{n}}{\beta_n} =1,
\end{align}
which is equivalent to saying 
\begin{equation}
s \leq {\gamma n}/{\beta_n}, \quad \forall s \leq |\hat{\Gamma}| \gamma/\hat{u}_n,  
\label{eq:sIsSmallForhatG}
\end{equation}
as desired.

For a set $S\subset \hat{\Gamma}$, we {now} characterize  the size of its neighborhood in $\hat{G}$, $|\mathcal{N}^*(S)|$, which depends on the distribution of the random subset, $\Delta$. Fix {some} $s\in \N$ {with} $s\leq  |\hat{\Gamma}|\gamma/\hat{u}_n$. 
{From} Eq.~\eqref{eq:sIsSmallForhatG}, we know that $s \leq {\gamma n}/{\beta_n}$.
{Consider some $S\subset \hat\Gamma$ with $|S|=s$.
Using the expansion property of $g_n$, we have that $|\mathcal{N}(S)|\geq\beta_n s$. Therefore,} 
 \begin{align}
\pb(|\mathcal{N}^*(S)| \leq \hat{u}_n s) = & \pb\left(\sum_{j \in \neig{S}} \mathbb{I}(j \in \Delta) \leq \hat{u}_n s\right) \nln
\stackrel{(a)}{\leq} & \pb\left(\mbox{Bino}\lt(|\neig{S}|\, , \frac{b_n}{n} \left(\rho+3\epsilon/4\right) \rt) \leq \hat{u}_n s\right) \nln
\stackrel{(b)}{\leq} & \pb\left(\mbox{Bino}\lt(\beta_n s \, , \frac{b_n}{n} \left(\rho+3\epsilon/4\right) \rt) \leq \hat{u}_n s\right),
\label{eq:NstarS_bnd}
\end{align}
for all sufficiently large $n$. Step $(a)$ follows from {the assumption that $\pb\left(j \in \Delta\right)  {\geq} (\rho + 3\epsilon/4)\frac{b_n}{n}$}, and step $(b)$ from the inequality $|\mathcal{N}(S)|\geq\beta_n s$. We observe that
\begin{align}
\mu \bydef & \E\left(\mbox{Bino}\lt(\beta_n s \, , \frac{b_n}{n} \left(\rho+3\epsilon/4\right) \rt)\rt) \nnb\\
=& (\rho+3\epsilon/4) \frac{\beta_n b_n}{n} s\nnb\\
\stackrel{(a)}{=} & (\rho+3\epsilon/4) \frac{1}{n}\cdot {\frac{80}{\epsilon^2}}\cdot \frac{n \ln n}{\beta_n} \beta_n s \nnb\\
= & (\rho+3\epsilon/4)\frac{{80}\ln n}{{\epsilon^2}} s, 
\label{eq:expandmu}
\end{align}
where in step $(a)$ we used the substitution {$b_n = \frac{80}{\epsilon^2}\cdot \frac{n \ln n}{\beta_n}$}. We also have that
\begin{align}
\hat{u}_n  =& \beta_n  \frac{m_n}{n} \nnb\\
=& \beta_n  \frac{\rho b_n}{{\hat \rho} n} \nnb\\
=&  \beta_n \frac{\rho }{{\hat \rho} n} \cdot {\frac{80}{\epsilon^2}}\cdot \frac{n \ln n}{\beta_n}\nnb\\
= & \frac{\rho}{{\hat \rho}} \cdot \frac{{80}\ln n}{{\epsilon^2}}. 
\label{eq:expandalphas}
\end{align}

By combining Eqs.~\eqref{eq:expandmu} and \eqref{eq:expandalphas}, we can derive a useful lower bound {on} the quantity $1-\frac{s \hat{u}_n }{\mu}$, which is recorded in the lemma that follows.

\begin{lemma} 
\label{lem:chernoffGap1}
We have that
\begin{equation}
1- \frac{s \hat{u}_n }{\mu} \geq \frac{\epsilon}{2}. 
\label{eq:mualpharatio}
\end{equation}
\end{lemma}

\def\xx{3\epsilon/4}
\bpf 
Using Eqs.~\eqref{eq:expandmu} and \eqref{eq:expandalphas} in the first step below, we have that 
$$
1- \frac{s {\hat{u}_n} }{\mu}= 1-\frac{\rho}{\hat\rho(\rho+\xx)}.$$
Recall that $\epsilon=(1-\rho)/2$, so that $\rho=1-2\epsilon$ and that $\hat\rho=1/(1+\epsilon/4)$. Using these substitutions, we obtain
\begin{align*}
1- \frac{s \hat{u}_n }{\mu}=& 1-\frac{(1-2\epsilon) (1+\epsilon/4)}{1-2\epsilon+3\epsilon/4}\\
=&\frac{3\epsilon/4 -\epsilon/4+2\epsilon^2/4}{1-5\epsilon/4}\\
=&\frac{\epsilon(1+\epsilon)/2}{1-5\epsilon/4}\\
\geq&\frac{\epsilon}{2}.
\end{align*}

\qed

To obtain an upper bound for the probability in Eq.~\eqref{eq:NstarS_bnd}, we substitute Eqs.~\eqref{eq:expandmu} and \eqref{eq:mualpharatio} into Eq.~\eqref{eq:NstarS_bnd}.  {Given the assumption that} $s {\leq} {\gamma n}/{\beta_n}$, {we have that}
\begin{align}
\pb(|\mathcal{N}^*(S)| \leq \hat{u}_n s) \leq& \pb\left(\mbox{Bino}\lt(\beta_n s \, , \frac{b_n}{n} \left(\rho+3\epsilon/4\right) \rt) \leq \hat{u}_n s\right)  \nnb\\
\stackrel{(a)}{\leq} & \exp\left(-\frac{1}{2}\left(\frac{\epsilon}{2}\right)^2\mu\right) \nln
\stackrel{(b)}{=} & \exp\left(- \frac{\epsilon^2}{8} \cdot \frac{{80}\ln n}{{\epsilon^2}}   (\rho+3\epsilon/4)s\right) \nln
= & {\exp(-(10\ln n)(\rho+3\epsilon/4)s)} \nln
{\sk{c}{\leq}} & \exp(-(5\ln n) s)\nln
=& \frac{1}{n^{5s}}.
\label{eq:NstarS_goodbnd} 
\end{align}
for all sufficiently large $n$.  Step $(a)$ is based on a  multiplicative form of the Chernoff bound (cf.~Chapter 4 of \cite{mitzenmacher2005probability}),  $\pb\left(X \leq (1- \delta)\mu\right)\leq \exp\left(-\frac{1}{2}\delta^2\mu\right)$, where $X$ is a binomial random variable with $\EX{X}=\mu$, and 
\begin{equation}
\delta  = 1- \frac{s \hat{u}_n }{\mu} \geq \epsilon/2,
\end{equation}
where the {last} inequality follows from Lemma \ref{lem:chernoffGap1}. Step $(b)$ follows from Eq.~\eqref{eq:expandmu}, and $(c)$ from the assumption that $\rho\geq 1/2$.

We now apply Eq.~\eqref{eq:NstarS_goodbnd} to subsets of $\hat{\Gamma}$, and {use the union bound.} We have, for all sufficiently large $n$, that 
\begin{align}
 \pb \lt( \, \mbox{$\hat{G}$ is not a $\lt({\gamma}/{\hat{u}_n} \, , \hat{u}_n\rt)$-expander}  \, \rt)  {\leq} & \pb( \, \mbox{$\exists S \subset \hat{\Gamma}\mbox{ such that: } |S|\leq |\hat{\Gamma}|{\gamma}/\hat{u}_n$ {and} $|\mathcal{N}^*(S)|\leq \hat{u}_n|S|$}  ) \nln 
  \stackrel{(a)}{\leq}  & \sum_{s=1}^{|\hat{\Gamma}|\gamma/\hat{u}_n} \left(\sum_{S \subset \hat{\Gamma}, |S|=s}  \pb\left(|\mathcal{N}^*(S)|\leq \hat{u}_n s\right) \right) \nln
\leq & \sum_{s=1}^{|\hat{\Gamma}|\gamma/\hat{u}_n}  {|\hat{\Gamma}|\choose s} \pb\left(|\mathcal{N}^*(S)|\leq \hat{u}_n s\right) \nln 
  \stackrel{(b)}{<} & \sum_{s=1}^{|\hat{\Gamma}|\gamma/\hat{u}_n}  b_n^s \pb\left(|\mathcal{N}^*(S)|\leq \hat{u}_n s\right) \nln 
\stackrel{(c)}{\leq} & \sum_{s=1}^{|\hat{\Gamma}|\gamma/\hat{u}_n} b_n^s {\frac{1}{n^{5s}}} \nln
\leq& \sum_{s=1}^\infty {(b_n/n^5)^s} \nln
= & {\frac{b_n/n^5}{1-b_n/n^5}. }
\label{eq:probGhatExpand}
\end{align}
Step $(a)$ {is the} union bound. In step $(b)$, we used the bound  ${n \choose k} \leq n^k$, and the fact that $|\hat{\Gamma}| = m_n =\frac{\rho}{\hat \rho} b_n <b_n$. Step $(c)$ follows from Eq.~\eqref{eq:NstarS_goodbnd}. 
 Because {$\beta_n \gg \ln n$, we have that $b_n\lsim \frac{n\ln n}{\beta_n} \ll n$, and hence}
 \begin{equation}
\frac{b_n}{n^5}\leq \frac{1}{n^{3}}, 
\label{eq:theta_nsmall}
 \end{equation}
 for all sufficiently large $n$. Combining Eqs.~\eqref{eq:probGhatExpand} and \eqref{eq:theta_nsmall}, we conclude that 
 \begin{equation}
  \pb \lt( \, \mbox{$\hat{G}$ is not a $\lt(\frac{\gamma}{\hat{u}_n}, \hat{u}_n\rt)$-expander}  \, \rt) \leq \frac{1}{n^3}, 
 \end{equation}
 for all sufficiently large $n$. This proves our claim. \qed

\section{Expanded Modular Architectures}
\label{app:expanModular}
{In this appendix,}  we start by describing the graph product, and subsequently we discuss the implications of using an expander graph.

\bfpara{Construction of the Architecture.} We first express the average degree as a product, $d_n = d^m_n\cdot d^e_n$, where the relative magnitudes of $d^m_n$ and $d^e_n$ are a design choice. The architecture is constructed as follows.
\begin{enumerate}
\item Similar to the case of the Modular architecture, partition $I$ and $J$ into equal-sized clusters of size $d^m_n$. We will refer to the index set of the queue and server clusters as $\calQ$ and $\calS$, respectively. For any $i\in I$ and $j\in J$, denote by $q(i) \in \calQ$ and $s(j)\in \calS$,  the indices of the queue and server clusters to which $i$ and $j$ belong, respectively. 
\item Let $g^e_n$ be a bipartite graph of maximum degree $d^e_n$ 
whose left and right nodes are the
queue and server clusters, $\calQ$ and $\calS$, respectively. Let $E^e$ be the set of edges of $g^e_n$.
\item To construct the interconnection topology $g_n=(I\cup J,E)$, let $(i, j) \in E$ if and only if their corresponding queue and server clusters are connected in $g_n^e$, i.e., if $(q(i),s(j))\in E^e$. 
\end{enumerate}

Note that by the above construction, each queue is connected to at most $d^e_n$ server clusters through $g_n^e$, and within each connected cluster, to  $d^m_n$ servers. Therefore, the maximum degree of $g_n$ is $d^m_n\cdot d^e_n=d_n$.

\bfpara{Scheduling Policy.} The scheduling policy requires the knowledge of the arrival rate vector, $\blam_n$, and involves two stages. For a given $\blam_n$, the computation in the first stage is performed  only once, while the steps in the second stage are repeated throughout the operation of the system. 
\begin{enumerate}
\item Compute a feasible flow, $\{f_{q,s}\}_{(q,s)\in E^e}$, over the graph $g^e_n$, where the incoming flow at each queue cluster $q\in \calQ$ is equal to $\sum_{i\in q}\lambda_i$, and the outgoing flow at each server cluster $s\in \calS$ is constrained to be less than or equal to $\frac{1+\rho}{2} d^m_n$. ({It turns out  that, under our assumptions, such a feasible flow exists \cite{xu2014power}}.)  Denote by $f_{q ,s}$ the total rate of flow from the queue cluster $q$ to the server cluster $s$. 
\item {Arriving} jobs first wait in queue until they are fetched by a server. When a server becomes available, it chooses a neighboring queue cluster (w.r.t.\ the topology of $g_n^e$) with probability roughly proportional to the flow between the clusters. In particular, a server in cluster $s$  chooses the queue cluster $q$ with probability 
\begin{equation}
p_{s,q} = \frac{f_{q,s}}{\sum_{q' \in \neig{s}}f_{q',s}}\cdot \frac{1+\rho}{2}+\frac{1}{\deg(s)}\cdot \frac{1-\rho}{2},
\label{eq:psqdef}
\end{equation} 
where $\deg(s)$ is the degree of $s$ in $g_n^e$. Within the chosen cluster, the server starts serving a job from an arbitrary non-empty queue, or, if all queues in the cluster are empty, the server initiates an idling period whose length  is exponentially distributed with mean $1$. 
\end{enumerate}

When the graph $g^e_n$ is an {expander graph}, we refer to the topology created via the above procedure as an \emph{Expanded Modular architecture generated by $g_n^e$}. 

Note that an Expanded Modular architecture is constructed as a ``product'' between an expander graph across the queue and server clusters, and a fully connected graph for each pair of connected clusters. As a result, its performance is also of a hybrid nature: the expansion properties of $g^e_n$ guarantee a large capacity region, while a diminishing delay is obtained as a result of the growing size of the server and queue clusters. We summarize this in {the} next theorem. Here we assume that $d^e_n$ is sufficiently large so that the expander graph described in Lemma \ref{cor:expan} exists. The reader is referred to Section 3.4.5 of \cite{xu2014power} for {the} proof of the theorem {(although with different choices for some of the constants)}.

\begin{theorem}[Capacity and Delay of Expanded Modular Architectures]\label{thm:expandModular}
Suppose that $d_n=d^m_n\cdot d^e_n$. Let $\gamma =\sqrt{\rho}$ and $\beta_n = \frac{1}{2}\cdot  \frac{\ln (1/\rho)}{1+\ln (1/\rho)} \ d^e_n$. Let $g_n^e$ be a $(\gamma/\beta_n, \beta_n)$-expander with maximum degree $d^e_n$, and let $g_n$ be an Expanded Modular architecture generated by $g_n^e$. If
\begin{equation}
u_n\leq  \frac{1+\rho}{2}\beta_n = \frac{1+\rho}{4}\cdot  \frac{\ln (1/\rho)}{1+\ln (1/\rho)} \ d^e_n,
\end{equation}
then, under the scheduling policy described above, we have that
\begin{equation}
\sup_{\blam_n \in \bLam_n(u_n)} \EX{W | \blam_n} \lsim \frac{c}{ d^m_n},
\label{eq:expModuDim}
\end{equation}
where $c$ is a constant that does not depend on $n$. 
\end{theorem}


\emph{A Tradeoff between the Size of the Capacity Region and the Delay}. For the Expanded Modular architecture, the relative values of $d^m_n$ and $d^e_n$ reflect a design choice: a larger value of $d^e_n$ ensures a larger capacity region, while a larger value of $d^m_n$ yields smaller delays. Therefore, while the Expanded Modular architecture is able to provide a strong delay guarantee that applies to \emph{all} arrival rate vectors in $\bLam_n(u_n)$, it comes at the expense of either a slower rate of diminishing delay (small $d^m_n$) or a smaller capacity region (small $d^e_n$).

\end{APPENDICES}

\end{document}

%% file: macroDef.tex
\newcommand{\beq}{\begin{equation}}
\newcommand{\eeq}{\end{equation}}
\newcommand{\beqn}{\begin{eqnarray}}
\newcommand{\eeqn}{\end{eqnarray}}
\newcommand{\benum}{\begin{enumerate}}
\newcommand{\eenum}{\end{enumeratse}}

\newcommand{\etmz}{\end{itemize}}

\newcommand{\bthm}{\begin{thm}}
\newcommand{\ethm}{\end{thm}}
\newcommand{\bdefn}{\begin{defn}}
\newcommand{\edefn}{\end{defn}}
\newcommand{\lt}{\left}
\newcommand{\rt}{\right}
\newcommand{\nnb}{\nonumber}

\newcommand{\bydef}{\stackrel{\triangle}{=}}

\newcommand{\smid}{\, | \,}

\global\long\def\bfpara#1{{\bf #1} }

\global\long\def\bpf{\emph{Proof.}$\,$ }

\global\long\def\epf{\qed}

\newcommand{\mcal}{\mathcal}

\global\long\def\del#1{}

\global\long\def\bydef{\stackrel{\triangle}{=}}

\global\long\def\lsim{\lesssim}
\global\long\def\gsim{\gtrsim}

\global\long\def\sk#1#2{\stackrel{(#1)}{#2}}

\global\long\def\calQ{\mathcal{Q}}
\global\long\def\calS{\mathcal{S}}

\newcommand{\pb}{\mathbb{P}}
\newcommand{\N}{\mathbb{N}}
\newcommand{\R}{\mathbb{R}}

\newcommand{\zp}{{\mathbb{Z}_+}}
\global\long\def\rp{\mathbb{R}_{+}}

\global\long\def\llfl{}
\global\long\def\rrfl{}

\global\long\def\EX#1{\mathbb{E}\left(#1\right)}
\global\long\def\E{\mathbb{E}}
\global\long\def\var#1{\mbox{\normalfont Var}\left(#1\right)}

\global\long\def\neig#1{\mathcal{N}\left( #1 \right)}

\global\long\def\bLam{\mathbf{\Lambda}}

\global\long\def\gcal{\mathcal{G}}

\global\long\def\blam{\bm{\lambda}}

\global\long\def\lun{u_n}
\global\long\def\ldn{d_n}

\global\long\def\degr{\mbox{\normalfont deg}}

\global\long\def\bigo#1{\mathcal{O}\left({#1}\right)}
\newcommand{\BR}{\mathbf{R}}

\global\long\def\km1{t_{k-1}}

\global\long\def\iden{\mathbb{I}}

\global\long\def\l1#1{\left\|#1\right\|_\infty}


%% file: VirtualQueueSection_V15_Rev2-Final.tex
\section{Analysis of the Expander Architecture}
\label{sec:VirtualQueues}

In this section,  we introduce a policy for the Expander architecture, based on batching and virtual queues, which will then be used  to prove Theorem \ref{thm:expanderArch}. We begin by describing the basic idea at a high level. 

\subsection{The Main Idea}
\label{sec:virtual}
Our policy proceeds by collecting a fair number of arriving jobs to form batches. Batches are thought of as being stored in a virtual queue,  with each batch treated as a single entity. By choosing the batch size large enough, one expects to see certain statistical regularities that can be exploited in order to efficiently handle the jobs within a batch. We now provide an outline of the operation of the policy, for a special case.

Let us fix $n$ and consider the case where $\lambda_{i}=\lambda<1$ for all $i$. Suppose that at time $t$, all servers are busy serving some job. Let us also fix some $\gamma_n$ such that $\gamma_n\ll 1$, while $n\gamma_n$ is large. During the time interval $[t,t+\gamma_n)$, ``roughly'' $\lambda n \gamma_n$ new jobs will arrive and $n\gamma_n$ servers will become available. Let  $\Gamma$ be the set of queues that received any job and let $\Delta$ be the set of  servers that became available during this interval. Since $\lambda n\gamma_n \ll n$, these incoming jobs are likely to be spread out across different queues, so that most queues receive at most one job. Assuming that this is indeed the case, we focus on  $g_n|_{\Gamma\cup \Delta}$, that is, the connectivity graph $g_n$, restricted to $\Gamma \cup \Delta$. The key observation is that this is a subgraph sampled uniformly at random among  all subgraphs of $g_n$ with approximately $\lambda n \gamma_n$ left nodes and $n\gamma_n$ right nodes. When $n \gamma_n$ is sufficiently large, and $g_n$ is well connected (as {in an {expander} with appropriate expansion properties}), we expect that, with high probability, $g_n|_{\Gamma\cup \Delta}$ admits a  matching that includes the entire set $\Gamma$ (i.e., a one-to-one mapping from $\Gamma$ to $\Delta$). In this case,  we can ensure that \emph{all} of the roughly $\lambda n \gamma_n$ jobs can start receiving service at the end of the interval, by assigning them to the available servers in $\Delta$ according to this particular matching. Note that the resulting queueing delay will be comparable to $\gamma_n$, which has been assumed to be small.

The above described scenario corresponds to the normal course of events. However, with a small probability, the above scenario may not materialize,
due to statistical fluctuations, such as:
 
\begin{enumerate}
\vspace{-3pt}
\item Arrivals may be concentrated on a small number of queues.
\vspace{-3pt}
\item The servers that become available may be located in a subset of $g_n$ {that is not well connected to the queues with arrivals.} 
\end{enumerate}
In such cases, it may be impossible to assign the jobs in $\Gamma$ to servers in $\Delta$. These exceptional cases will be handled by the policy in a different manner. However, if we can guarantee that the probability of such cases is low, we can then argue that their impact on performance is negligible.

Whether or not the above mentioned exceptions will have low probability of occurring depends on whether the underlying connectivity graph, $g_n$, has the following property: with high probability, a randomly sampled sublinear (but still sufficiently large) subgraph of $g_n$ admits a large set of ``flows.'' This property will be used to guarantee that, with high probability, the jobs in $\Gamma$ can {indeed} be assigned to distinct servers in the set~$\Delta$. {We will show that an {expander} graph with appropriate expansion does possess this property.}

\subsection{An additional assumption}
Before proceeding, we introduce an additional assumption on the arrival rates, which will remain in effect throughout this section, and which will simplify some of the arguments. Appendix \ref{app:LowBonArrival} explains why {this} assumption can be made without loss of generality. 

\begin{assumption}
\label{ass:lamLowerBound}
{\bf (Lower Bound on the Total Arrival Rate)} We have that $\rho\in (1/2 \, , \, 1)$, and the total arrival rate satisfies the lower bound
\begin{equation}
\sum_{i =1}^n \lambda_i\geq (1-\rho)n. 
\label{eq:sumLamLowerBound}
\end{equation}
\end{assumption}

\subsection{The Policy}
\label{sec:thePolicy}

We now describe in detail the scheduling policy. Besides $n$, the scheduling policy uses the following inputs: 
\begin{enumerate}
\item  $\rho$, the traffic intensity introduced in Condition \ref{cond:arrdis}, in Section \ref{sec:ArchiPerMetri},
\item $\epsilon$, a positive constant such that $\rho+\epsilon<1$.
\item $b_n$, a batch size parameter,
\item $g_n$, the connectivity graph. 
\end{enumerate}
Notice that the {arrival rates, $\lambda_i$, and the} fluctuation parameter, $u_n$, {are} \emph{not}  {inputs} to the scheduling policy. 

At this point it is useful to make a clarification regarding the $\lsim$ notation. Recall that the relation $f(n)\lsim g(n)$ means that $f(n)\leq c g(n)$, for all $n$, where $c$ is a positive constant.
Whenever we use this notation, we require that the constant $c$ cannot depend on any parameters other than $\rho$ and $\epsilon$. Because we view $\rho$ and $\epsilon$ as fixed throughout, this makes $c$ an absolute constant. 

\subsubsection{Arrivals of Batches.} \label{s:barr}
Arriving jobs are organized in \emph{batches} of cardinality $\rho\hspace{1pt} b_n$, where $b_n$ is a design parameter, to be specified later.\footnote{In a slight departure from the earlier informal description, we define batches by keeping track of the number of arriving jobs as opposed to keeping track of time.}
Let $T^B_0=0$. For $k\geq 1$, let $T^B_k$ be the time of the $(k\rho b_n)$th arrival to the system, which we also view as the {\it arrival time} of the $k$th batch. For $k\geq 1$,  the $k$th batch consists of the $\rho b_n$ jobs that arrive during the time interval $(T^B_{k-1},T^B_k]$. The length $A_k=T^B_k-T^B_{k-1}$ of this interval will be called the {\it $k$th inter-arrival time}. We record, in the next lemma, some immediate statistical properties of the batch inter-arrival times.
 
\begin{lemma}
\label{lem:arr}
The batch inter-arrival times,  $\left\{A_k\right\}_{k\geq 1}$, are i.i.d., with $$
\frac{b_n}{n} \leq 	\E\left(A_k\right)\leq \frac{\rho}{1-\rho}\cdot \frac {b_n}{n},$$ 
and $\var{A_k} \lsim b_n/n^2$. 
\end{lemma}
\bpf  
The batch inter-arrival times are i.i.d., due to  our independence assumptions on the job arrivals. By definition, $A_k$ is equal in distribution to the  time until a Poisson process  records $\rho b_n$ {arrivals}. 
This Poisson process has rate  $r=\sum_{i=1}^n\lambda _i$, and using also Assumption  \ref{ass:lamLowerBound} in the first inequality below,
we have
$$(1-\rho)n \leq \sum_{i=1}^n\lambda_i  {=r}\leq \rho n,$$
The random variables $A_k$ are Erlang  (sum of $\rho b_n$ exponentials {with rate $r$). Therefore,} 
$$\EX{A_k} {=\rho b_n \cdot \frac{1}{r}} \geq \rho b_n\cdot \frac{1}{\rho  n} = \frac{b_n}{n}.$$
Similarly,
$$\E(A_k){=\rho b_n \cdot \frac{1}{r}} \leq \rho b_n \cdot \frac{1}{(1-\rho) n}.$$
Finally,
$$\var{A_k}  {=\rho b_n \cdot \frac{1}{r^2}}
\leq \llfl \rho b_n\rrfl \cdot \frac{1}{(1-\rho)^2 \llfl  n\rrfl^2} \lsim \frac{b_n}{ n^2}.$$ \epf

\subsubsection{The Virtual Queue}
\label{s:virtual}
Upon arrival, batches are placed in what we refer to as a \emph{virtual queue}.  The virtual queue is a GI/G/1 queue, which is operated in FIFO fashion. That is, a batch waits in queue until all previous batches are served, and then starts being served by {a} virtual queueing system. The service of a batch by the virtual queueing system lasts until a certain time {by} which all jobs in the batch have already been assigned to, and have started receiving service from, one of the physical servers, at which point 
the service of the batch is completed and the batch departs  from the virtual queue. The time elapsed from the initiation of service of batch until its departure is called the \emph{service time} of the batch.
As a consequence, the queueing delay of a job in the actual (physical) system is bounded above by the sum of:
\begin{itemize}
\item[(a)] the time from the arrival of the job until the arrival time of the batch that the job belongs to; 
\item[(b)] the time that the batch  waits in the virtual queue;
\item[(c)] the service time of the batch.
\end{itemize}

\noindent
\paragraph{\bf Service slots.} The service of the batches at the virtual queue is organized along consecutive time intervals that we refer to as \emph{service slots.} The service slots are intervals of the form $(\ell\s, (\ell+1)\s]$, where $\ell$ is a nonnegative integer, whose length is\footnote{To see how the length of the service slot was chosen, recall that the size of each batch is equal to $\rho \hspace{1pt}b_n$. The length of the service slot hence ensures that $(\rho+\epsilon)b_n$, the expected number of servers that will become available (and can therefore be assigned to jobs) during a single service slot, is  greater than the size of a batch, so that there is hope of assigning all of these jobs to available servers within a single service slot. At the same time, since $\rho+\epsilon<1$, service slots are shorter than the expected {batch} inter-arrival time, which is needed  for the stability of the virtual queue.}
$$\s=(\rho+\epsilon)\cdot\frac{b_n}{n}.$$

We will arrange matters so that batches can complete service and depart \emph{only} at the end of a service slot, that is, at times of the form $\ell\s$. Furthermore, we assume that the physical servers are operated as follows. If  either a batch completes service at time $\ell\s$ or if there are no batches present at the virtual queue at that time, we assign to every idle server a dummy job whose duration is an independent exponential random variable, with mean 1. This ensures that the state of the $n$ servers is the same (all of them are busy) at certain special times, thus facilitating further analysis, albeit at the cost of some inefficiency.

\subsubsection{The Service Time of a Batch.} 
The specification of the service time of a batch depends on whether the batch, upon arrival, finds an empty or nonempty virtual queue.

\begin{figure}[h]
\centering
\includegraphics[scale=.52]{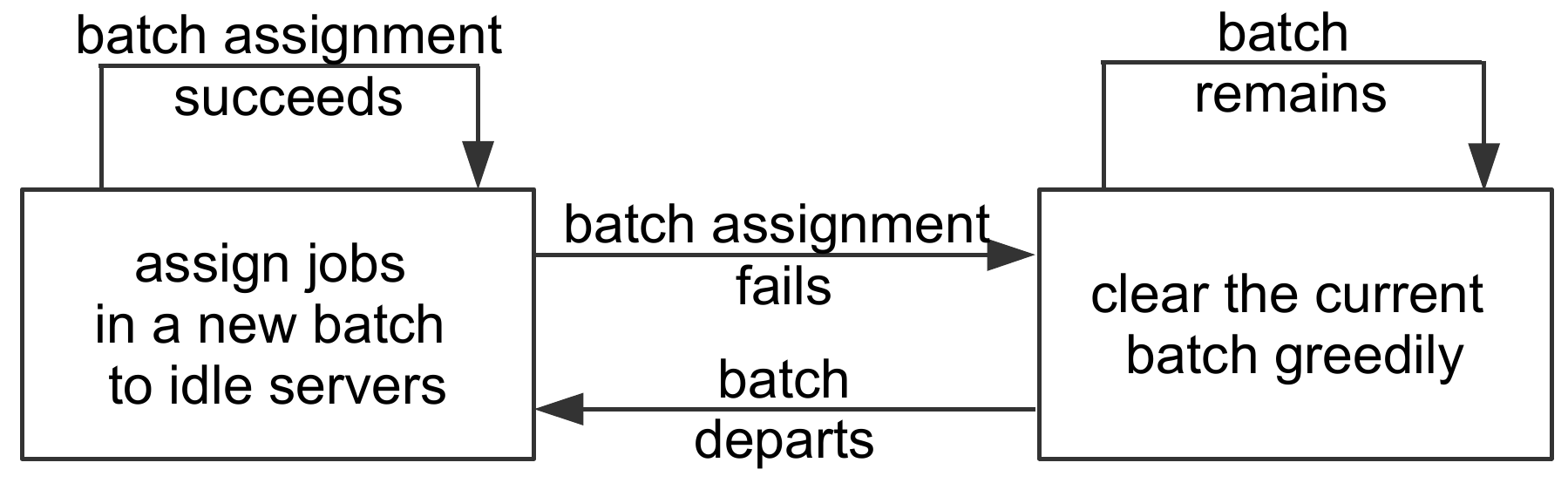}
\caption{{An illustration of the service slot dynamics. An arrow indicates the transition from the end of one service slot to the next.}}
\label{fig:serviceSlots}
\vspace{-15pt}
\end{figure}

Suppose that a batch arrives during the service slot $(\ell\s, (\ell+1)\s]$ and finds an empty virtual queue; that is, all previous batches have departed by time $\ell \s$. 
According to what was mentioned earlier, at time $\ell\s$, all physical servers are busy, serving either real or dummy jobs. 
Up until the end of the service slot, any server that completes service is not assigned a new (real or dummy) job, and remains idle, available to be assigned a job at the very end of the service slot. Let $\Delta$ be the set of servers that are idle at time $(\ell+1)\s$, the end of the service slot. At that time, we focus on the jobs in the batch under consideration. We wish to assign each job $i$ in this batch to a {\it distinct} server $j\in \Delta$, subject to the constraint that $(i,j)\in E$. We shall refer to such a job-to-server assignment as a \emph{batch assignment}. There are two possibilities {(cf.~Figure \ref{fig:serviceSlots})}: 
\begin{itemize}
\item[(a)]
If a batch assignment can be found, each job in the batch is assigned to a server according to that assignment, and the batch departs at time $(\ell+1)\s$. In this case, we say that the service time of the batch was \emph{short}.
\item[(b)] If a batch assignment cannot be found, 
we start assigning the jobs in the batch to physical servers in some arbitrary greedy manner: Whenever a server $j$ becomes available, we assign to it a job from the batch under consideration, and from  some queue $i$ with $(i,j)\in E$, as long {as} such a job exists. (Ties are broken arbitrarily.)
As long as every queue is connected to at least one server,  all jobs in the associated batch will be eventually assigned. The 
last of the jobs in the batch gets assigned during a subsequent service interval $(\ell'\s,(\ell'+1)\s]$, where $\ell'>\ell$, and we define $(\ell'+1)\s$ as the departure time of the batch.
\end{itemize}
If the $k$th batch did indeed find an empty virtual queue upon arrival, its service time, denoted by $S_k$, is  the time elapsed from its arrival until its departure.

Suppose now that a batch arrives during a service slot $(\ell\s, (\ell+1)\s]$ and finds a non-empty virtual queue; that is, there are one or more batches that arrived earlier and which have not departed by time $\ell\s$. In this case, the batch waits in the virtual queue until some time of the form $\ell'\s$, with $\ell'>\ell$, when the last of the previous batches departs. 
Recall that, as specified earlier,  at time $\ell'\s$ all servers are made to be busy (perhaps, by giving them dummy jobs) and we are faced with a situation identical to the one considered in the previous case, as if  
the batch under consideration just arrived at time $\ell'\s$; in particular, the same service policy can be applied. For this case, where the $k$th  batch arrives to find a non-empty virtual queue, its service time, $S_k$, extends from the time of the departure of the $(k-1)$st batch until the departure of the $k$th batch.

\subsection{Bounding the Virtual Queue by a GI/GI/1 Queue}
\label{s:modif}
Having defined the inter-arrival {and service} times of the batches, the virtual queue is a fully specified,  work-conserving, FIFO single-server queueing system.

We note however one complication. The service times of the different batches are dependent on the arrival times. To see this, suppose, for example, that a batch upon arrival sees an empty virtual queue and that its service time is ``short.'' {Then,} its service time will be equal to the remaining time  until the end of the current service slot, and therefore dependent on the batch's arrival time.
Furthermore, the service times of different batches are dependent: if the service time of the previous batch happens to be too long, then the next batch is likely to see upon arrival a non-empty virtual queue, which then implies that its own service time will be an integer multiple of $\s$.

In order to get around these complications, and to be able to use results on GI/GI/1 queues, we define the \emph{modified service time,} $S'_k$, of the $k$th service batch to be equal to $S_k$, rounded above to the nearest integer multiple of $\s$:
$$S_k'=\min\{\ell\s : \ell\s\geq S_k,\ \ell=1,2,\ldots\}.$$
Clearly, we have $S_k\leq S_k'$.

We now consider a modified (but again FIFO and work-conserving) virtual queueing system in which  the arrival times are the same as before, but the service times are the $S_k'$. A simple coupling argument, based on Lindley's recursion, shows that for every sample path, the time that the batch spends waiting in the queue of the original virtual queueing system  is less than or equal to the time spent waiting in the queue of the modified virtual queueing system. It therefore suffices to upper bound the expected time spent in the queue of the modified virtual queueing system. 

We now argue that the modified virtual queueing system is  a GI/GI/1 queue, i.e., that the service times $S'_k$ are i.i.d., and independent from the arrival process. For a batch whose service  starts during the {service} slot $[\ell\s,(\ell+1)\s)$,  the modified service time is equal to $\s$, whenever the batch service time is short. Whether the batch service time will be short or not is determined by the composition of the jobs in this batch and by the identities of the servers who complete service during the service slot $[\ell\s,(\ell+1)\s)$.
Because the servers start at the same ``state'' (all busy) at each service slot, it follows  that the events that determine whether a batch service time will be short or not are independent across batches, and with the same associated probabilities. 

Similarly, if a batch service time is not short, the additional
time to serve the jobs in the batch is affected only by the composition of jobs in the batch and the service completions at the physical servers after time $\ell\s$, and these are again independent from the inter-arrival times and the modified service times $S_m'$ of other
batches $m$. Finally, the same considerations show the independence of the $S_k'$ from the {batch} arrival process.

It should now be clear from the above discussion that the modified service time of a  batch is of the form
\begin{equation} 
S_k'=\s+X_k\cdot \hat S_k,
\label{eq:Smdecomp}
\end{equation}
where:
\begin{itemize}
\item[(a)]
$X_k$ is a Bernoulli random variable which is equal to 1 if and only if the $k$th batch service time is not short, i.e., it takes more than a single service slot;
\item[(b)]
$\hat S_k$ is a random variable which (assuming that every queue is connected to at least one server) is stochastically dominated by the sum of $\rho b_n$ independent exponential random variables with mean 1, rounded up to the nearest multiple of $\s$. (This dominating random variable corresponds to the extreme case where all of the $\rho b_n$ jobs in the batch are to be served in sequence, by the same physical server.)
\item[(c)] The pairs $(X_k,\hat S_k)$ are i.i.d.
\end{itemize}

\subsection{Bounds on the Modified Service Times}

For the remainder of Section \ref{sec:VirtualQueues}, we will assume that \begin{equation}
\mbox{$g_n$ is a $(\gamma/\beta_n, \beta_n)$-expander,}
\end{equation}
where $\gamma$ and $\beta_n$ are defined as in the statement of Theorem \ref{thm:expanderArch}.  

The main idea behind the rest of the proof is as follows.
We will upper bound the expected time spent in the modified virtual queueing system using Kingman's bound \cite{kingman61} for GI/GI/1 queues. Indeed, the combination of a batching policy with Kingman's {bound} is a fairly standard technique for deriving delay upper bounds (see, e.g., \cite{TS08}). 
We already have bounds on the mean and variance of the inter-arrival times. In order to apply Kingman's bound, it remains to obtain bounds on the mean and variance of the service times $S'_k$ of the modified virtual queueing system. 

We now introduce an important quantity associated with a graph $g_n$, by defining 
$$q(g_n)=\pb\left(X_k = 1\smid g_n\right);$$
because of the i.i.d.\ properties of {the} batch service times, this quantity does not depend on $k$.
In words, for a given connectivity graph $g_n$, the quantity $q(g_n)$ stands for the probability that we cannot find a batch assignment, between the jobs in a batch and the servers that become idle during a period of length $\s$.

We begin with the following lemma, which provides bounds on the mean and variance of $S_k'$.

\begin{lemma} \label{lem:mv}
{There exists a sequence, $\{c_n\}_{n \in \N}$, with $c_n\lsim b_n$, such that for all $n\geq 1$}
$$ {s \leq} \EX{ S_k' \mid  g_n} \leq  \s+q(g_n)c_n,$$
and
$$\var{S_k' \mid  g_n}\lsim q(g_n) c_n^2.$$
\end{lemma}
\bpf
{The fact that $\EX{ S_k' \mid  g_n} \geq s$ follows from the definition of $S_k'$ in Eq.~\eqref{eq:Smdecomp} and the non-negativity of $X_k \hat{S}_k$.} The definition of an expander ensures that every queue is connected to at least one server through $g_n$. Recall that $\hat S_k$ is zero if $X_k=0$; on the other hand, if $X_k=1$, and as long as
every queue is connected to some server, {then} $\hat S_k$ 
is  upper bounded by the sum of $\rho b_n$ exponential random variables with mean 1, rounded up to an integer multiple of $\s$. 
Therefore, 
$${\EX{ S_k' \mid  g_n}  = s+\EX{X_k\hat S_k \mid g_n}
=s+ \pb(X_k = 1\smid g_n) \cdot\EX{\hat S_k \mid X_k=1,g_n}}
\leq s+q(g_n)(\rho b_n +\s),$$ which leads to the first bound in the statement of the lemma, with $c_n=b_n+\s$. Since $\s$ is proportional to $b_n/n$, we also have $c_n\lsim b_n$, as claimed. 
Furthermore,
$$\var{S_k' \mid g_n}=\var{X_k \hat S_k \mid g_n}
\leq \EX{{X_k^2} \hat S_k^2 \mid g_n} \lsim q(g_n) (b_n+\s)^2 =q(g_n) c_n^2.$$
\qed

We now 
need to obtain bounds on $q(g_n)$. This is nontrivial and forms the core of the proof of the theorem.  In what follows, we will show that with appropriate assumptions on the various parameters, and for any $\blam\in\bLam(u_n)$, an Erd\"os-R\'enyi random graph has a very small $q(g_n)$, with high probability.

\subsection{Assumptions on the Various Parameters}\label{s:params}

{From now on,} we focus on a specific batch size {parameter} of the form
\begin{equation}
b_n= \frac{{320}}{(1-\rho)^2}\cdot\frac{n\ln n}{\beta_n}.
\label{eq:b_n}
\end{equation}
We shall also set 
\begin{equation}
\epsilon = \frac{1-\rho }{2}. 
\label{eq:setEpsilon}
\end{equation}

We assume, as in the statement of Theorem \ref{thm:expanderArch}, that $d_n\ll n$, and {that} 
\begin{equation}
\beta_n \gsim d_n {\gg \ln n}. 
\label{eq:dnconsDef}
\end{equation}
Under these choices of $b_n$ and $d_n$,  we have
\begin{equation}
b_n\lsim \frac{n}{{d_n/ \ln n}}\ll n;
\label{eq:bnissmalln}
\end{equation}
that is, the batch size is \emph{vanishingly small} compared to $n$. Finally, we {will only consider arrival rate vectors that belong} to the set $\bLam_n(u_n)$ (cf.\ Condition \ref{cond:arrdis}), where, as in the statement of Theorem \ref{thm:expanderArch},
\begin{equation}
\lun \leq {\frac{1-\rho }{2} \beta_n}. 
\label{eq:unconsDef}
\end{equation}

\subsection{The Probability of a Short Batch Service Time}

We now come to the core of the proof, aiming to show that 
{if the connectivity graph $g_n$ is an expander graph with a sufficiently large expansion factor, then} 
$q(g_n)$ is small. More precisely, we aim to show that a typical batch will have high probability of having a short service time. A  concrete statement is given in the result that follows, and the rest of this subsection will be devoted to its proof.

\begin{proposition}
\label{prop:expanderSmallq}
Fix $n\geq 1$. We have that
\begin{equation}
q(g_n) \leq \frac{1}{n^2}. 
\end{equation}
\end{proposition}

\def\bA{{\bf A}}

Let us focus on a particular batch, and let us examine what it takes for its service time to be short. There are two sources of randomness: 
\begin{enumerate}
\item A total of $\rho b_n$ jobs arrive {to the queues}. Let $A_i$ be the number of jobs that arrive at the $i$th queue, let $\bA=(A_1,\ldots,A_n)$, and let $\Gamma$ be the set of queues that receive at least one job. In particular, we have
$$\sum_{i=1}^n A_i = \sum_{i\in \Gamma} A_i = \rho b_n.$$
\item During the time slot at which the service of the batch starts, each server starts busy (with a real or dummy job). With some probability, and independently from other servers or from the arrival process, a server becomes idle by the end of the service time slot. Let $\Delta$ be the set of servers that become idle.
\end{enumerate}
Recalling the definition of $X_k$ as the indicator random variable of the event that the service time of the $k$th batch is not short, we see that $X_k$ is completely determined by the graph {$g_n$} together with $\bA$ and $\Delta$. For the remainder of this subsection, we suppress the subscript $k$, since we are focusing on a particular batch. We therefore have a dependence of the form 
$$X=f(g_n, \bA, \Delta),$$
for some function $f$, and we emphasize the fact that {$\bA$} and {$\Delta$} are independent.

Recall that $\epsilon=(1-\rho)/2$, and from the statement of Theorem \ref{thm:expanderArch} that
\begin{equation}\label{eq:rr}
\hat{\rho}= \frac{1}{1+(1-\rho)/8} = \frac{1}{1+\epsilon/4}.
\end{equation}
Clearly, $\hat\rho <1$, and with some elementary algebra, it is not difficult to show that, for {any given} $\rho \in (0,1)$, 
$$\hat{\rho} > \rho.$$
Let 
$$m_n =\frac{\rho}{\hat{\rho}} b_n,$$
so that
\begin{equation}\label{eq:rmn}
\hat{\rho} m_n=\rho b_n.
\end{equation}
Finally, let
$$
\hat{u}_n= \beta_n \frac{m_n}{n}.$$

We will say that $\bA$ is \emph{nice} if there exists a set $\hat\Gamma \supset  \Gamma$ of cardinality $m_n$,  such that  $A_i=0$ whenever $i\notin\hat\Gamma$, and
$$A_i{<} \hat{u}_n,\qquad\forall\ i \in \hat{\Gamma}.$$

We now establish that $\bA$ is nice, with high probability. 
The main idea is simple: $\bA$ is not nice only if one out of a finite collection of binomial variables with large means  takes a value which is away from its mean by a certain multiplicative factor. Using  the Chernoff bound, this probability can be shown to decay {at least} as fast as $1/n^3$. The details of this argument are given in the proof of Lemma \ref{lem:nice}, in Appendix \ref{app:lem:nice}. 

\begin{lemma} \label{lem:nice}
For all sufficiently large $n$, we have that
$$\pb(\bA \mbox{\rm\ is not nice}) \leq \frac{{1}}{n^3}.$$
\end{lemma}

We now wish to establish that when $\bA$ is nice, there is high probability (with respect to $\Delta$), that the batch service time will be short.
Having a short batch service time is, by definition, equivalent to the existence of a batch assignment, which in turn is equivalent to the existence of a certain flow in a subgraph of $g_n$. The lemma that follows deals with the latter existence problem {for the original graph, but will be later applied to subgraphs.}  Let $\overline{\BR}(g)$ be the closure of the capacity region, $\BR(g)$, of $g$.
 
\begin{lemma} 
\label{lem:expandFlowClosure}
Fix $n, n' \in \N$, $\rho\in (0,1)$, and $\gamma >\rho$. Suppose that {an $n\times n'$ bipartite graph}, $g_n$, is a {$\left({\gamma}/{\beta_n}, \beta_n\right)$-expander,  where $\beta_n\geq u_n$}. Then  $\bLam_n(u_n) \subset \overline{\BR}(g_n)$. \end{lemma}

\bpf
The claim follows directly from Lemma \ref{lem:expandRobust}, by noting that $\overline{\BR}(g_n)\supset \BR(g_n)$. \qed 

The next lemma is the key technical result of this subsection. It states that if $g_n$ is an expander, then, for any given $\hat\Gamma$, the random subgraph $g_n|_{\hat\Gamma \cup \Delta}$ will be an expander graph with high probability (with respect to $\Delta$). {The lemma is stated as a stand-alone result, though we will use a notation that is consistent with the rest of the section.} The proof {relies on a delicate application of the Chernoff bound}, and is given in Appendix \ref{app:lem:subGraphExpansion}. 

\begin{lemma}
\label{lem:subGraphExpansion}
Fix $n\geq 1$, $\gamma \in (0,1)$, {and $\rho\in [1/2 \, , 1)$. Let $g_n = (I\cup J, E)$ be an $n\times n$ bipartite graph that is a $(\gamma/\beta_n, \beta_n)$-expander, where $\beta_n \gg \ln n$. Define the following quantities: 
\begin{align}
\epsilon =& \frac{1-\rho}{2}, \nln
\hat{\rho} =& \frac{1}{1+\epsilon/4}, \nln
b_n =  & \frac{320}{(1-\rho)^2}\cdot \frac{n\ln n}{\beta_n} = \frac{80}{\epsilon^2}\cdot \frac{n\ln n}{\beta_n}, \nln
m_n = & \frac{\rho}{\hat{\rho}} b_n, \nln
\hat{u}_n = &  \beta_n\frac{m_n}{n}.
\label{eq:expanLemmaQuant}
\end{align}}
{Let $\hat{\Gamma}$ be an arbitrary subset of the left vertices, $I$, such that
\begin{equation}
|\hat{\Gamma}|= m_n,
\end{equation} 
and let $\Delta$ be a random subset of the right vertices, $J$, where each vertex belongs to $\Delta$ independently and with the same probability, where
\begin{equation}
\pb(j \in \Delta) \geq (\rho+3\epsilon/4)\frac{b_n}{n}, \quad \forall j\in J, \label{eq:jdlt1}
\end{equation}
for all $n$ sufficiently large.  Denote by $\hat{G}$ the random subgraph $g_n|_{\hat{\Gamma} \cup \Delta}$.} Then
\begin{equation}
 \pb \lt( \, \mbox{$\hat{G}$ is not a $\lt(\gamma / \hat{u}_n \, , \hat{u}_n\rt)$-expander}  \, \rt) \leq \frac{1}{n^3}, 
\end{equation} 
for all $n$ sufficiently large, where the probability is measured with respect to the randomness in $\Delta$. 
\end{lemma}

\def\ba{{\bf a}}

{To invoke Lemma \ref{lem:subGraphExpansion}, note that the conditions in Eq.~\eqref{eq:expanLemmaQuant} are identical to the definitions for the corresponding quantities in this section. We next verify that Eq.~\eqref{eq:jdlt1} is satisfied by the random subset, $\Delta$, consisting of the idle servers at the end of a service slot. Recall that the length of a service slot is $\frac{b_n}{n} \left({\rho+\epsilon}\right)$, and hence the probability that a given server, $j$, becomes idle by the end of a service slot is
\begin{align}
\pb\left(j \in \Delta\right)  =  1-\exp\left(-\frac{b_n}{n} \left(\rho+\epsilon\right)\right) {\sim}  \left(\rho+\epsilon\right)\frac{b_n}{n},
\end{align}
as $n\to \infty$. Therefore, for all $n$ sufficiently large, we have that $\pb\left(j \in \Delta\right)  {\geq} (\rho + 3\epsilon/4)\frac{b_n}{n}. $}
We will now apply Lemmas \ref{lem:expandFlowClosure} and \ref{lem:subGraphExpansion} to the random subgraph with left (respectively, right) nodes $\hat\Gamma$ (respectively $\Delta$), and with the demands $A_i$, for $i\in \hat\Gamma$, playing the role of $\blam$.

\begin{lemma}
\label{lem:short}
If $n$ is large enough, and if the value {$\ba$} of $\bA$ is nice, then
$$\pb(X=1\mid \bA=\ba) \leq  \frac{1}{n^3},$$
where the probability is with respect to the randomness {in $\Delta$.}
\end{lemma}

\bpf
We fix some $\ba$, assumed to be nice. Recall that $$\hat{u}_n = \beta_n\frac{m_n}{n},$$ and from the statement of Theorem \ref{thm:expanderArch} that $$\gamma=\sqrt{\hat{\rho}} > \hat{\rho}.$$
We apply  Lemma \ref{lem:expandFlowClosure} to the {randomly sampled subgraph} $\hat G $, with left nodes $\hat\Gamma$, $|\hat\Gamma|=m_n$, and right nodes $\Delta$. We have the following correspondence: the parameters $n$ and $\rho$,  in 
Lemma \ref{lem:expandFlowClosure} become, in the current context, 
$m_n$  and $\hat\rho$, respectively, and the parameters $\beta_n$ and $u_n$ {both} become $\hat{u}_n$. Thus, by Lemma \ref{lem:expandFlowClosure},
\begin{equation}
\mbox{if $\hat{G}$ is a $(\gamma/\hat{u}_n, \hat{u}_n)$-expander, \ \ then\ \   }\bLam_{m_n}(\hat{u}_n) \subset \overline{\BR}(\hat{G}). 
\label{eq:blam_mnInBRGhat}
\end{equation}

Let $\hat\bA$ be the vector of job arrival numbers $\bA$, restricted to the set of nodes in $\hat\Gamma$, and let $\hat\ba$  be the realization of   $\hat\bA$. Note that we have
$$\sum_{i\in\hat\Gamma}{\hat a}_i=\sum_{i=1}^n {A_i}  \ =\rho b_n=\hat{\rho} m_n,$$
because of Eq.\ \eqref{eq:rmn}. Furthermore, for any $i\in\hat\Gamma$, the fact that  $\ba$ is nice implies that  $\hat a_i<\hat{u}_n$. 
{Thus,} $\hat\ba \in {\bLam_{m_n}(\hat{u}_n)}$. By Eq.~\eqref{eq:blam_mnInBRGhat}, this further implies that
\begin{equation}
\mbox{if $\hat{G}$ is a $(\gamma/\hat{u}_n, \hat{u}_n)$-expander, \ \ then\ \   }
\hat \ba \in \overline{\BR}(\hat{G}).\label{eq:ainbLammn}
\end{equation}

By Lemma \ref{lem:subGraphExpansion}, the graph $\hat{G}$ is a $(\gamma/\hat{u}_n, \hat{u}_n)$-expander with probability at least $1-n^{-3}$. Combining this fact with  Eq. \eqref{eq:ainbLammn}, we have thus verified that $\hat\ba$ belongs to $\overline{\BR}(\hat{G})$, with probability at least
$$1-\frac{1}{n^3}.$$  
With $\overline{\BR}(\hat{G})$ {having been}  defined {as} the closure of the capacity region, $\BR(\hat{G})$ (cf.~Definition \ref{def:deFlow}), the fact that the vector $\hat \ba$ belongs to $\overline{\BR}(\hat{G})$ is a statement about the existence of a {feasible} flow, $
\{\hat f_{ij}: (i,j)\in {\hat{E}}\}$ {(where $\hat E$ is the set of edges in $\hat G$)},
in a linear network flow model {of the form}
\begin{align}
\hat a_i = \sum_{j: (i,j)\in {\hat E}} f_{ij}, &\quad \forall i \ \in \hat \Gamma, \nln
\sum_{i: (i,j)\in \hat{E}} f_{ij} \leq 1, &\quad \forall\  j \in  {\Delta},  \nln
f_{ij} \geq 0,& \quad \forall \ (i,j)\in {\hat E }.  \nnb
\end{align}
Because the ``supplies'' $\hat a_i$ in this network flow model, as well as the unit capacities of the right nodes are integer, it is well known that there also exists an integer flow. That is, we can find  $\hat f_{ij}\in\{0,1\}$ such that $\sum_j \hat f_{ij}=\hat a_i$, for all $i$, and $\sum_i \hat f_{ij}\leq 1$, for all $j$. But this is the same as the statement that there exists a feasible batch assignment  {over $\hat G$}.  Thus, {for large enough $n$ and for any given  nice $\ba$, the conditional} probability that a batch assignment does not exist is upper bounded by ${n^{-3}}$, as claimed.
\qed

We can now complete the proof of Proposition 
\ref{prop:expanderSmallq}. By considering unconditional probabilities where $\bA$ is random, and for $n$ large enough,
we have that
\begin{align}\label{eq:x}
\pb(X=1)\leq & \pb(\bA \mbox{{\rm \ is not nice}})
+\sum_{\ba\ {\rm nice}} \pb(X=1\mid \bA=\ba)\cdot \pb(\bA=\ba) \nln
\sk{a}{\leq} & \frac{{1}}{{n^3}} +\sum_{\ba\ {\rm nice}} \pb(X=1\mid \bA=\ba) \cdot \pb(\bA=\ba) \nln
\sk{b}{\leq} & \frac{{1}}{{n^3}} +\sum_{\ba\ {\rm nice}} \frac{1}{{n^3}}\cdot \pb(\bA=\ba) \nln
\leq& \frac{{2}}{{n^3}} \nln
\leq& \frac{1}{n^2},
\end{align}
where steps $(a)$ and $(b)$ follow from Lemmas \ref{lem:nice} and \ref{lem:short}, respectively. This concludes the proof of Proposition \ref{prop:expanderSmallq}. 

%


\subsection{Service and Waiting Time Bounds for the  Virtual Queue}

\subsubsection{Service Time Bounds.}
We will now use Lemma \ref{lem:mv} and {Proposition \ref{prop:expanderSmallq}} to bound the mean and variance of the service times  in the modified virtual queue.

\begin{lemma} \label{lem:SMstat} The {modified batch} service times, $S'_k$, are i.i.d., with $$
	\EX{S'_k {\,|\, g_n}} \sim  (\rho+\epsilon)\cdot {\frac{b_n}{n}},\ \  {\rm and} \ \ \var{S'_k {\,|\, g_n}} \lsim {\frac{b_n^2}{n^2}}.$$ 
\end{lemma}

\bpf 
We use the fact from Lemma \ref{lem:mv}, that 
$s\leq \EX{S'_k {\,|\, g_n}}\leq s+q(g_n)c_n$, where $c_n\lsim b_n$.
We recall that $s=(\rho+\epsilon)b_n/n$,  and use the fact $q(g_n){\leq} {n^{-2}}$, as guaranteed by Proposition \ref{prop:expanderSmallq}. The term $q(g_n)c_n$ satisfies $q(g_n)c_n\lsim b_n/{n^2}$, which is of lower order than $b_n/n$, and hence negligible compared to $s$. This proves the first part of the lemma.

For the second part, we use Lemma \ref{lem:mv} in the first inequality below, and the fact that $q(g_n){\leq} {n^{-2}}$ in the second, to obtain
$$\var{S'_k {\,|\, g_n}}\lsim q(g_n) c_n^2 \lsim{\frac{b_n^2}{n^2}}.
$$
\qed

\subsubsection{Waiting Time Bounds.}
{Fix $n$ and the graph $g_n$. Let}
$W^B$ be a random variable whose distribution is the same as the steady-state distribution of the time that a batch spends waiting in the queue of the  virtual queueing system introduced in Section \ref{s:virtual}.

\begin{prop} \label{prop:mdelay} 
We have that
\begin{equation}
	\EX{{W^B} { \,|\, g_n}} \lsim {\frac{b_n}{n}}. 
\end{equation}
\end{prop}

\bpf
As discussed in Section \ref{s:modif}, {the
waiting time of a batch, in the virtual queueing system, is dominated by the waiting time in a modified virtual queueing system, which} is a GI/GI/1 queue, with independent inter-arrival times $A_k$  (defined in Section \ref{s:barr}) and independent service times $S'_k$. 
{Let $W'$ be a random variable whose distribution is the same as the steady-state distribution of the time that a batch spends waiting in the queue of the  modified virtual queueing system.}

According to Kingman's bound 
\cite{kingman61}, $W'$ satisfies
$$\EX{W'\,|\, g_n}\leq \tilde{\lambda} \frac{\sigma_a^2 + \sigma_s^2}{2(1-\tilde{\rho})},$$ where $\tilde{\lambda}$ is the arrival rate, $\tilde{\rho}$ is the traffic intensity, and $\sigma_a^2$ and $\sigma_s^2$ are the variances of the inter-arrival times and service times, respectively, that are associated with the modified virtual queueing system.

From Lemma \ref{lem:arr}, we have
$$\tilde{\lambda} = \frac{1}{\EX{A_k}} \leq \frac{n}{b_n}. $$
and
$$\sigma_a^2 = \var{A_k}\lsim  \frac{b_n}{n^2}.$$
We now bound
$$\tilde{\rho} = \frac{\EX{{S'_k\,|\, g_n}}}{\EX{A_k}}.
$$
From the first part of Lemma \ref{lem:SMstat}, we have $\EX{S'_k\,|\, g_n}\sim(\rho+\epsilon)b_n/n$. Together with the bound $1/\EX{A_k} \leq n/b_n$, we obtain that as 
$n\to \infty$, $\tilde{\rho}$ is upper bounded by a number strictly less than 1. We also have, from the second part of Lemma \ref{lem:SMstat},
$$\sigma_s^2 = \var{S'_k {\mid g_n}} \lsim  \frac{b_n^2}{n^2}.$$
Using these inequalities in Kingman's bound, we obtain
 {$$
\EX{W^B\mid g_n}\leq\EX{W' \mid g_n} \lsim { \frac{n}{b_n} \cdot \frac{b_n^2}{n^2}=
 \frac{b_n}{n}.}$$
 }
 \qed

\subsection{Completing the Proof of Theorem \ref{thm:expanderArch}}
\label{sec:thm1proof}
\bpf  As discussed in Section \ref{s:virtual}, the expected waiting time of a job is upper bounded by the sum of three quantities.
\begin{itemize}
\item[(a)] The expected time from the arrival of the job until the arrival time of the batch that the job belongs to. This is bounded above by the expected time until there are $\rho b_n$ subsequent arrivals, which is equal to $\EX{A_1}$. By Lemma \ref{lem:arr}, this is bounded above  by  $c_1 b_n/n$,  for some constant $c_1$. 
\item[(b)] The expected time that the batch waits in the virtual queue. This is also upper bounded by $c_2 b_n/n$, by Proposition \ref{prop:mdelay},  for some constant $c_2$. 
\item[(c)] The service time of the batch, which (by Lemma \ref{lem:SMstat}) again admits an upper bound of the form $c_3 b_n/n$,  for some constant $c_3$. 
\end{itemize}
Furthermore, in the results that give these upper bounds,  $c_1$, $c_2$, and $c_3$, are absolute constants, that do not depend on $\blam_n$ or $g_n$. 

By our assumptions on the choice of $b_n$ in Section \ref{s:params}, we have
$b_n= \frac{{320}}{(1-\rho)^2}\cdot\frac{n\ln n}{\beta_n}$,
{and $\beta_n$ is proportional to $d_n$.} We conclude that there exists a constant $c$ such that for large enough $n$, we have
$\E\left(W\,|\,g_n,\blam_n\right) \leq  c \ln n/d_n$, for any given $\blam_n\in\bLam_n(u_n)$, which is an upper bound of the desired form. This establishes {P}art $1$ of the theorem. Finally, {P}art $2$ follows from the way that the policy was constructed.
\epf

\subsection{On Practical Policies} 
\label{sec:simu}

\begin{figure}
\centering
\includegraphics[scale=.3]{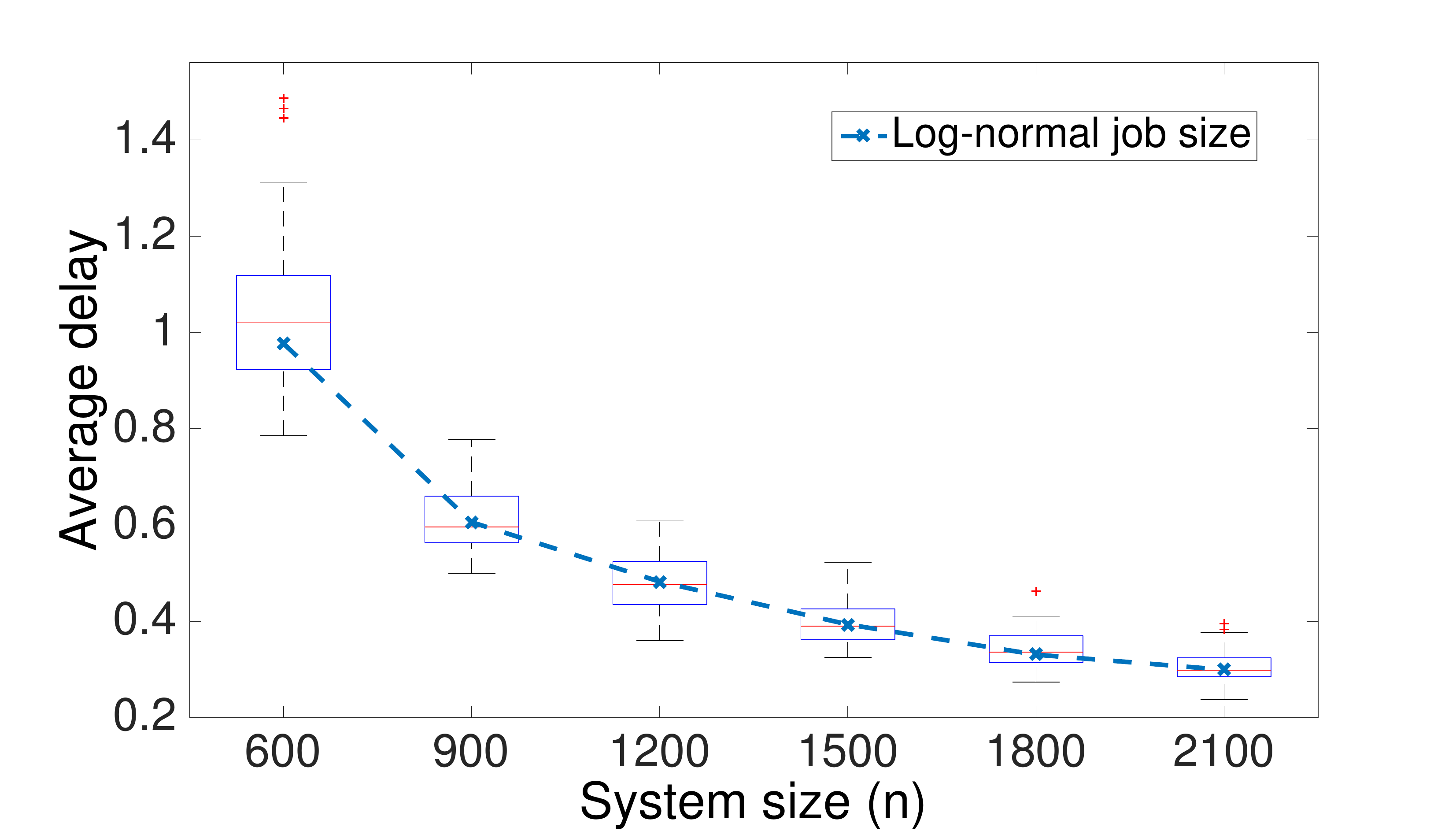}
\caption{{Simulations of the virtual-queue based policy given in Section \ref{sec:thePolicy}, with $d_n =  n^{2/3}$, $b_n = n\ln(n)/d_n$, and $\lambda_i = 0.5$ for all $i=1,\ldots, n$. The boxplot contains the average delay from 50 runs of simulations where the job size distribution is assumed to be exponential with mean $1$. Each run is performed on a random $d_n$-regular graph over $10^4$ service slots and a $1000$-slot burn-in period. The center line of a box represents the median and upper and lower edges of the box represent the 25th and 75th percentiles, respectively. The dashed line depicts the median average waiting times when the job sizes are distributed according to a log-normal distribution with mean $1$ and variance $10$.}}
\label{fig:simulations}
\vspace{-10pt}
\end{figure}

{Figure \ref{fig:simulations} provides simulation results for the average delay under the virtual-queue based scheduling policy used in proving Theorem \ref{thm:expanderArch}. The main role of the policy  is to demonstrate the fundamental potential of the Expander architecture in jointly achieving a small delay and large capacity region when the system size is large. In smaller systems, however, there could be other policies that yield better performance. For instance, simulations suggest that a seemingly naive greedy heuristic can achieve a smaller delay in moderately-sized systems, which is practically zero in the range of parameters in Figure \ref{fig:simulations}. Under the greedy heuristic, an available server simply fetches a job from a longest connected queue, and a job is immediately sent to a connected idle server upon arrival if possible. Intuitively, the greedy policy can provide a better delay because it avoids the overhead of holding jobs in queues while forming a batch. Unfortunately, it appears challenging to establish rigorous delay or capacity guarantees for the greedy heuristic and other similar policies. } 

{In some applications, such as call centers, the service times or job sizes may not be exponentially distributed (\cite{brown2005statistical}).  In Figure \ref{fig:simulations}, we also include the scenario where the job sizes are drawn from a log-normal distribution (\cite{brown2005statistical}) with an increased variance. Interestingly, the average delay appears to be somewhat insensitive to the change in job size distribution. }

